\documentclass[2pt]{article}
\usepackage{graphicx} 
\usepackage{subcaption}
\usepackage{fancyhdr} 
\usepackage{color} 
\usepackage{amsmath}
\usepackage{amsthm}
\usepackage{amssymb}
\usepackage{mathptmx}
\usepackage{appendix}
\usepackage{float}

\usepackage[margin=2cm]{geometry}
\usepackage{xcolor}
\usepackage{natbib}
\usepackage[linesnumbered,ruled,vlined]{algorithm2e}
\usepackage{longtable}
\usepackage{array}
\usepackage{booktabs}
\usepackage{caption}
\newcolumntype{P}[1]{>{\raggedright\arraybackslash}p{#1}}
\usepackage{authblk}

\author[1]{Ning Xie}
\author[1]{Lei Wei}
\author[2]{Nikola Bešinović}
\author[1]{Meng Wang}

\affil[1]{Chair of Traffic Process Automation, Technische Universität Dresden, Dresden 01069, Germany}
\affil[2]{Chair of Railway Operations, Technische Universität Dresden, Dresden 01069, Germany}

\title{Link-Based Multimodal Traffic Dynamics Model in Continuous-Time Framework}
\begin{document}
\fancyhf{}
\renewcommand{\headrulewidth}{0pt}
\cfoot{\thepage}

\maketitle
\begin{abstract}
Computationally efficient models for multimodal traffic flows with inter-modal interactions are foundational for coordinated traffic management. However, such models are lacking in the current literature. This study introduces a multimodal link transmission model (M-LTM) accommodating both continuous road traffic and discrete tramway traffic at the network level. M-LTM builds on a link model that captures the inter-modal interactions through the moving bottleneck theory. We further develop node models describing the flow transfer, tram operations, and mode interactions at tram stops and intersections. Specifically, the tram dwell process and dwell-induced congestion of road traffic are simulated at the stop node, and the intersection node model can reproduce queue spillback and tram priority.  Case studies are conducted on two synthetic one-node networks, a synthetic arterial network, and the real network of Dresden to demonstrate the predictive power and the applicability of the proposed model in network traffic analysis and management. The GEH statistics of road traffic are below 4, and the average tram arrival time absolute errors are less than 10 seconds.
\end{abstract}

\section{Introduction}
\label{introduction}

Modeling traffic evolution across networks is indispensable for urban traffic assignment, management, and control \citep{wu1998continuous,ahmed2018optimum}. The increasing prevalence of heterogeneous transport modes and complex interactions between them on shared infrastructures necessitates the development of multimodal traffic flow models \citep{martinez2023value}. Owing to the distinct operational characteristics, road traffic and tramway traffic have traditionally been modeled using segregated approaches. While effective within their respective domains, such disjoint modeling paradigms pose challenges for integrated multimodal traffic flow analysis and control \citep{fu2020empirical,patwary2021metamodel}. To address this, the study develops a unified traffic flow model that captures multimodal traffic dynamics and inter-modal interactions within a consistent continuous-time framework.

\subsection{Related work}
\label{related work}
Road traffic dynamics have been modeled in both microscopic and macroscopic approaches. Compared to microscopic ones \citep{nagel1992cellular,treiber2000congested}, macroscopic models \citep{dotoli2006urban,patwary2021metamodel,kang2022fractional,molnar2024destroying,mo2024game} enable aggregate representation of traffic flow and efficient computation for large-scale networks. Among these, the Lighthill–Whitham–Richards (LWR) model \citep{lighthill1955kinematic,richards1956shock} is one of the most influential macroscopic traffic flow models due to its solid theoretical foundation and ability to capture shockwave propagation.

\cite{lebacque1996godunov} shows the equivalence between first-order macroscopic traffic flow models and Godunov’s scheme, introducing a unified framework for modeling boundary conditions in the LWR model. This work lays the foundation for solving the first-order hyperbolic partial differential equations of LWR, leading to simplified formulations such as the Cell Transmission Model (CTM) and Link Transmission Model (LTM). The CTM discretizes space into cells and applies conservation laws to approximate traffic dynamics \citep{daganzo1994cell,daganzo1995cell}. However, its accuracy is limited in highly congested conditions unless a fine spatial discretization is used, which increases computational burden \citep{carey2014implementing}. In contrast, the LTM \citep{yperman2005link,yperman2006multicommodity,yperman2007link,raadsen2019continuous} represents traffic flow as continuum and describes dynamics across entire links based on Newell's simplified kinematic wave theory \citep{newell1993simplified}, allowing more accurate and computationally efficient modeling of shockwave propagation by avoiding the numerical diffusion introduced by coarse spatial discretization. 

Multiple adaptations of the LTM have been proposed to enhance its applicability in different traffic scenarios and modeling settings. Computational accuracy and efficiency have been improved through iterative solution algorithms \citep{himpe2016efficient}, while the assumption of triangular fundamental diagrams has been relaxed to general concave functions \citep{van2017extending,raadsen2016efficient}. Other extensions incorporate important traffic features, including variable speed \citep{hajiahmadi2013variable}, mixed traffic involving connected and automated vehicles \citep{zhang2023stochastic,lu2024link}, and turn-level queue representation \citep{wei2025link}. Another important research direction focuses on continuous-time formulations of the LTM \citep{jin2015continuous,han2016continuous,tumash2022multi,bliemer2019continuous}, which provide greater analytical tractability and can improve the accuracy of traffic flow representation. The continuous-time formulation is also adopted in this study.

The node model is another critical component of the LTM, describing traffic transfer at intersections \citep{smits2015family}. \cite{tampere2011generic} and \cite{jabari2016node} establish fundamental invariance and holding-free principles for node models and propose solutions for both signalized and unsignalized intersections. \cite{corthout2012non} further considers traffic conflicts at intersections and discusses solution uniqueness. For signalized intersections, substantial efforts have been devoted to capturing operational characteristics such as spillback conditions \citep{gibb2011model,han2015continuum}, phase switching \citep{durlin2008dynamic,canudas2019average}, and phase interactions \citep{yahyamozdarani2023continuous}. However, these studies mainly focus on private road traffic and generally model spillback effects only through downstream receiving constraints on connecting incoming links. As a result, they overlook the blockage of other conflicting movements caused by residual vehicles occupying the intersection area. Moreover, existing LTM-based node models are not well suited for public transport operations, particularly tram systems operating on dedicated tracks with strong timetable adherence. Important operational processes, such as tram dwelling and dwelling-induced blocking effects, are not explicitly captured in existing formulations.

Tramway traffic operates on partially dedicated Right-Of-Way in urban networks \citep{vuchic2002urban}. Along most links, trams travel on exclusive tracks with relatively stable operational patterns, while interactions with other road traffic mainly occur at stops, intersections, and other shared infrastructures. Unlike road traffic that can often be approximated as continuum, tramway traffic is dominated by discrete vehicle operations with relatively large headways. Consequently, tramway traffic dynamics are commonly modeled using microscopic approaches \citep{cuniasse2015analyzing}. Event-based models represent tram operations as sequences of events triggered by changes in vehicle states \citep{corman2014review,cacchiani2014overview,burghout2006discrete,toledo2010mesoscopic,sharma2023review}. These models provide foundations for tram scheduling and control optimization using disjunctive graphs \citep{d2009advanced,corman2010centralized,shakibayifar2018simulation}. Car-following models, such as Cellular Automata (CA) model \citep{li2005cellular,xun2013impact,qian2022bidirectional}, microscopically model the tramway traffic dynamics through operational interactions \citep{quaglietta2020multi,saidi2023train}, and have been extended with stochastic factors for operational variability and disturbances \citep{weik2020quantifying,corman2021stochastic}. However, these microscopic approaches are computationally demanding for large-scale multimodal traffic networks and are less suitable for representing aggregate tramway traffic dynamics and link/node-level flow propagation.

To evaluate network performance and simplify computational complexity, macroscopic approaches have been developed for railway or multimodal traffic \citep{liu2015integrating,gentile2017formulation}. In particular, the Macroscopic Fundamental Diagram (MFD) has been extended to railway traffic, including the Line Fundamental Diagram proposed for the Paris metro system \citep{cuniasse2015analyzing} and 3D-MFDs incorporating railway traffic flow, network traffic density, and passenger flow to describe delay and congestion. Further, the railway MFDs have been complemented by consideration of uncertainty \citep{weik2022macroscopic}, microscopic simulation \citep{dong2023description}, and Petri net \citep{goverde2007railway, szymula2024quantifying} for applications such as capacity analysis. The MFD-based approaches also capture aggregate input and output interactions between modes across the entire network \citep{loder2017empirics,balzer2023dynamic}. However, due to their aggregate features, these macroscopic approaches generally lack the ability to represent link/node-level operations and to explicitly model queue formation and propagation. Several studies extended CTM and LTM to public traffic by considering the urban traffic as mixed traffic flow \citep{liu2015integrating,bayrak2021optimization}, while they primarily focused on bus and can not be adapted to the discrete tramway flow because of the partially dedicated Right-Of-Way (ROW). Additionally, tram operations and interactions with road traffic are more complex. The tram dwelling processes are tightly scheduled following the First In First Out (FIFO) rule, and block all road traffic at unprotected stops. These operations and interactions were also overlooked in the CTM or LTM-based models for buses.

In summary, despite extensive developments in both road traffic LTM formulations and tramway operation models, a unified macroscopic framework capable of representing multimodal traffic dynamics and inter-modal interactions is still lacking. Existing LTM studies primarily focus on road traffic and do not explicitly account for tramway operational processes, such as dwelling at stops. In addition, current node models fall short in representing tram signal priority and queue spillback, limiting their applicability for large-scale multimodal traffic analysis and coordinated traffic control.

\subsection{Contribution}
\label{contribution}
To address the limitations of the current traffic dynamic models, this study extends the multimodal traffic dynamics analysis based on kinematic wave theory \citep{xie2026continuum} to the network level and develops a multimodal link transmission model (M-LTM). The model simulates road and tramway traffic flow evolutions with their interactions over networks within a unified continuous-time framework. Traffic propagation along links is depicted by the link model, while the node model simulates traffic transfer and interactions between the two traffic modes at tram stops and signalized intersections. Case studies are conducted on two one-node networks, an arterial network, and the Dresden network to verify the performance of the proposed model.

The proposed modeling framework advances the state-of-the-art in two ways.
    First, multimodal node models are developed to represent traffic transfer, tram operations, and inter-modal interactions for network-level traffic evolution. The proposed model explicitly captures tram dwelling operations, dwelling-induced road traffic blockage, and signalized intersection spillback. The node models enable a unified continuous-time multimodal LTM framework to represent traffic propagation and transfer over networks.
    Second, the proposed model is verified on synthetic multimodal traffic networks through comparisons with microscopic simulation results from SUMO, and on real-world multimodal networks with empirical data. Qualitative and quantitative evaluation metrics demonstrate that the proposed M-LTM can accurately reproduce multimodal traffic dynamics and inter-modal interactions with GEH statistics of road traffic below 4 and average tram arrival time absolute errors below 10 seconds.

The remainder of the paper is organized as follows. Section \ref{model overview} introduces the framework and the link model of the multimodal link transmission model, followed by the node model in Section \ref{node model} describing traffic dynamics between links. Several case studies are illustrated in Section \ref{case study} to verify the proposed model. Finally, Section \ref{conclusion} summarizes the study and demonstrates the future work.

\section{Model overview}
\label{model overview}
This section provides an overview of the proposed multimodal link transmission model in Section \ref{overall_framework}, followed by a revisit of the traffic propagation mechanism along links in Section \ref{link model}.
\subsection{Overall modeling framework}
\label{overall_framework}
\begin{figure}[!h]
    \centering
    \includegraphics[width=1\linewidth]{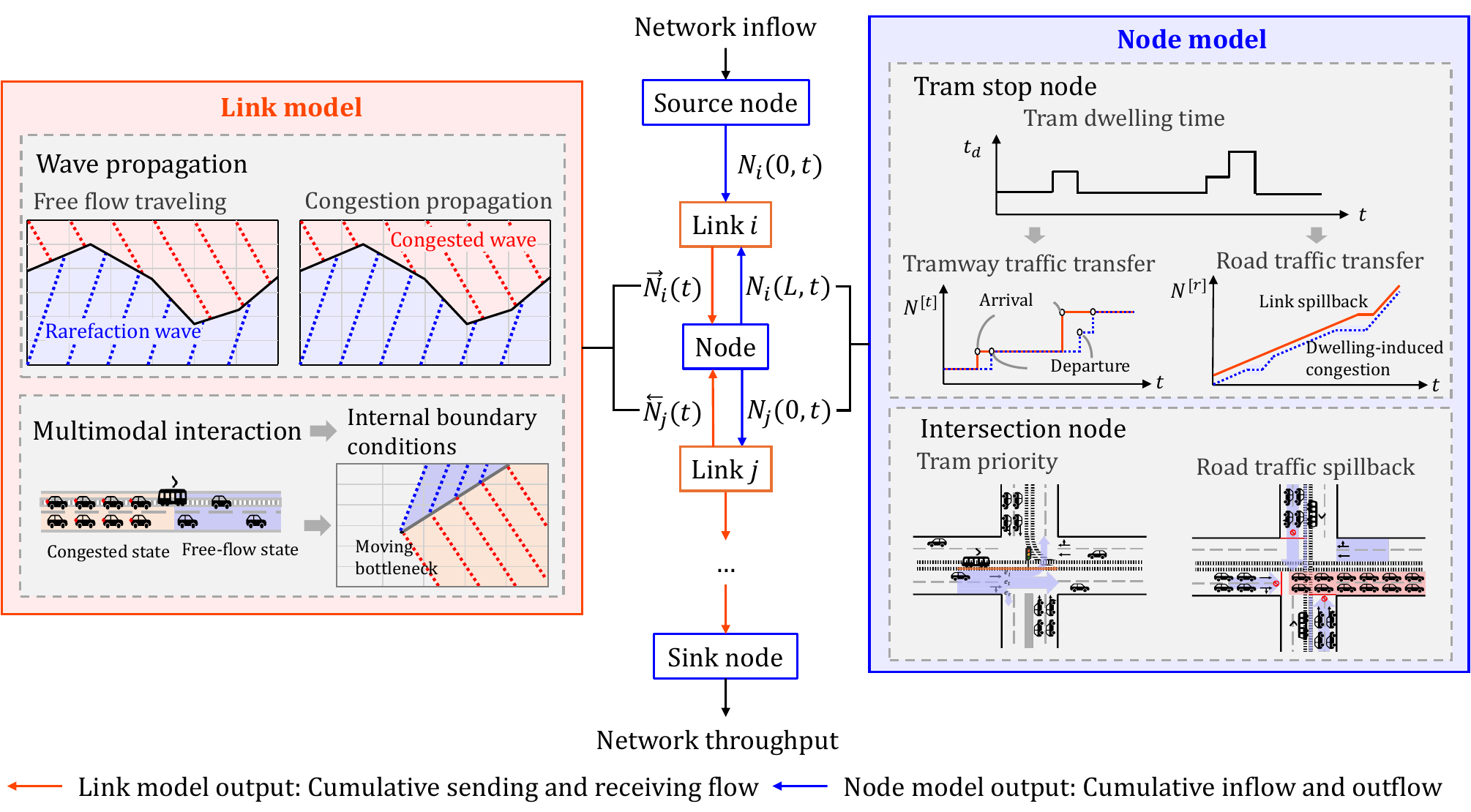}
    \caption{Framework of the Multimodal Link Transmission Model, where $N_i(0,t)$ and $N_j(0,t)$ are the cumulative inflows of link $i$ and $j$, $N_i(L,t)$ is the cumulative outflow of link $i$, and $\overrightarrow{N}_i(t)$ and $\overleftarrow{N}_j(t)$ are the cumulative sending flow of link $i$ and cumulative receiving flow of link $j$.}
    \label{framework}
\end{figure}

The proposed multimodal link transmission model simulates the road traffic and tramway traffic evolution in networks by the cumulative flows, comprising the link model and node model. The link model captures the rarefaction and congestion wave propagation along links, deriving the cumulative sending and receiving flows. The maximal cumulative flow that can leave the link at the downstream is referred to as cumulative sending flow. The maximal cumulative flow that can enter the link, assuming sufficient demand from the upstream node, is referred to as cumulative receiving flow. Additionally, the interactions between the two traffic modes on shared ROW links are incorporated using the moving bottleneck theory. The node model, classified into tram stop and signalized intersection nodes, achieves traffic transfer and determines the actual cumulative flow at the end of links. Specifically, the tram dwelling process is simulated at the stop, which impedes road traffic passage during dwelling time. The signalized intersection node model assigns traffic according to the demands from incoming links and supplies of outgoing links, and detects the spillback effects when the outgoing links are congested. Since the signal control is compulsory for tram safety at intersections, only signalized intersections are considered in this study. 
The model framework is illustrated in Figure \ref{framework}, clarifying the relationships between links and nodes, model variables, and model mechanisms. Links update the cumulative sending and receiving flows through the link model based on the cumulative inflows and outflows from source nodes or internal nodes, while the internal nodes update the cumulative inflow and outflow according to the cumulative sending and receiving flows of their upstream and downstream links. In the following text, the variables of road traffic and tramway traffic are distinguished by superscript ${[r]}$ and ${[t]}$, respectively.
\subsection{Traffic propagation on links}
\label{link model}
The link model is introduced in this section, which depicts the multimodal traffic evolution along homogeneous links based on the multimodal kinematic wave theory \citep{xie2026continuum}. We discuss two conditions where road and tramway traffic operate independently and the two traffic modes share the ROW in Section \ref{independent link model} and Section \ref{link model with shared ROW}. For simplicity, the link index subscript is omitted in this section without loss of clarity.

\subsubsection{General link model}
\label{independent link model}
The link model derives the cumulative sending and receiving flow at the downstream and upstream of the link according to the rarefaction and congestion wave propagation.

\begin{figure}[!h]
    \centering
    \includegraphics[width=0.7\linewidth]{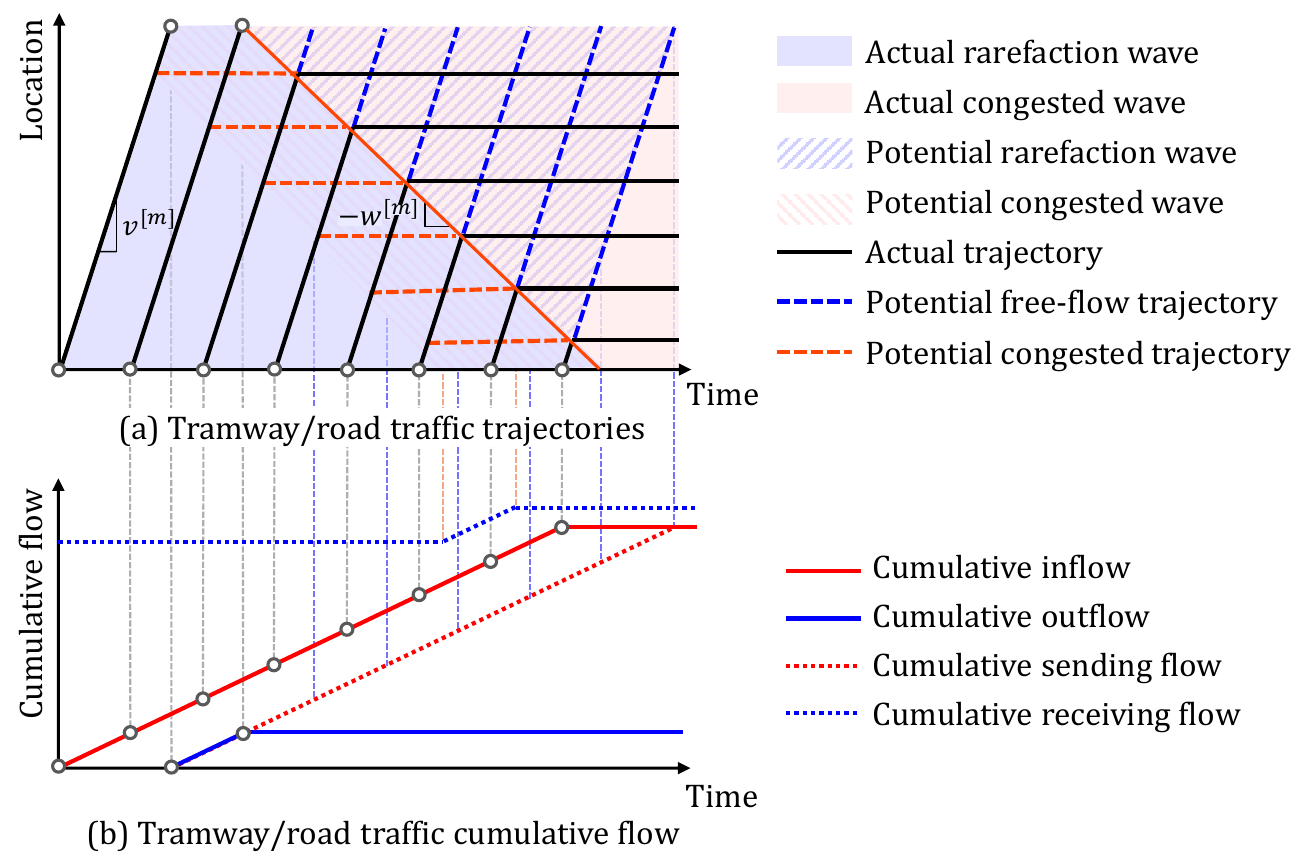}
    \caption{Illustration of the link model.}
    \label{link_model_illustration}
\end{figure}

In the general link model, road and tramway traffic operate on dedicated ROWs, without interactions on the link. Therefore, rarefaction and congestion wave propagation determine the traffic evolution along links. The cumulative sending flow is derived based on the cumulative inflow of the link, describing the rarefaction wave reaching the link downstream. Equation \ref{cumulative sending flow} formulates the cumulative sending flow, with the constraint of the saturation flow. The blue trapezoid in Figure \ref{link_model_illustration} (a) demonstrates the rarefaction wave propagation in the free-flow state, deriving the blue dashed curves of cumulative sending flow in Figure \ref{link_model_illustration} (b).
\begin{equation}
    \overrightarrow{N}^{[m]}(t)=\left\{
    \begin{array}{cc}
         0 & t-\frac{L}{v^{[m]}}<0 \\
         \min\{N^{[m]}(0,t-\frac{L}{v^{[m]}}),N^{[m]}(L,t-\frac{L}{v^{[m]}})+s^{[m]}\frac{L}{v^{[m]}}\}& t-\frac{L}{v^{[m]}}\geq 0
    \end{array}\right.
    \label{cumulative sending flow}
\end{equation}
\noindent
where $\overrightarrow{N}^{[m]}(t)$ is the cumulative sending flow of traffic mode $m$ (road traffic or tramway traffic) at time $t$, $N^{[m]}(0,t)$ is the cumulative inflow of traffic mode $m$ at time $t$, $L$ is the length of the link, $v^{[m]}$ is the free flow speed of $m$, and $s^{[m]}$ is the saturation flow of traffic mode $m$. The tramway traffic saturation flow is determined by the number of tracks and the minimum time gap.

Similarly, the shockwave from the link downstream propagates to the upstream, leading to congestion on the link as depicted by the red trapezoid in Figure \ref{link_model_illustration} (a). The cumulative receiving flow, indicated by red dashed curves in Figure \ref{link_model_illustration} (b), is calculated based on the cumulative outflow as Equation \ref{cumulative receiving flow}.
\begin{equation}
    \overleftarrow{N}^{[m]}(t)=\left\{
    \begin{array}{cc}
         \rho_j^{[m]}L & t-\frac{L}{w^{[m]}}<0 \\
         \min\{N^{[m]}(L,t-\frac{L}{w^{[m]}})+\rho_j^{[m]}L,N^{[m]}(L,t-\frac{L}{w^{[m]}})+s^{[m]}\frac{L}{w^{[m]}}\}& t-\frac{L}{w^{[m]}}\geq 0
    \end{array}\right.
    \label{cumulative receiving flow}
\end{equation}
\noindent
where $\overleftarrow{N}^{[m]}(t)$ is the cumulative receiving flow of traffic mode $m$ at time $t$, $N^{[m]}(L,t)$ is the cumulative outflow of traffic mode $m$ at time $t$, $w^{[m]}$ is the congestion wave speed of $m$, and $\rho_j^{[m]}$ is the congested density of $m$. The congested density of tramway traffic is given based on the relationship between effective occupancy and density, in which the tram effective occupancy is denoted by $\ell^{[t]}$ with the unit of $m/veh$.
Details of the multimodal traffic kinematic wave theory is referred to \cite{xie2026continuum}.

\subsubsection{Link model with inter-modal interactions}
\label{link model with shared ROW}

When the road space is limited, the two traffic modes may share the ROW. Trams at a lower speed reduce the link capacity for road traffic, which is regarded as moving bottlenecks \citep{newell1998moving,daganzo2005moving,daganzo2005numerical,munoz2002moving}. The moving bottleneck not only influences the road traffic dynamics, but is also influenced by the surrounding road traffic under congestion \citep{leclercq2004moving,simoni2017fast,wu2019estimating,li2023network}. These interactions are modeled as upstream and downstream internal boundary conditions in the link model, denoted by $\bar{N}_u^{[m]}$ and $\bar{N}_d^{[m]}$. If the tramway traffic free flow speed is the same as that of road traffic, the road traffic operates without influence, while the internal boundary conditions for tramway traffic are considered.

The tram is activated as a moving bottleneck when it slows down the incoming traffic from upstream. This happens when the road traffic demand is higher than the reduced capacity in the tram moving coordinate. Under such conditions, road traffic passes the tram in the relative flow $q_r$ in the tram moving frame. Conversely, when the road traffic flow is below the relative flow in the tram moving frame, the moving bottleneck is inactive. Specifically, if the congested road traffic travels at a lower speed than tramway traffic speed, the tram is blocked by surrounding traffic. Considering these interactions, we formulate the moving bottleneck influence as internal boundary conditions in three steps. First, the tram entry time $t_b$ is identified and used to determine whether the congested wave caused by the tram reaches the link upstream at time $t$. This allows us to compute the road traffic downstream internal boundary condition. Then, the tram influence on the link downstream is checked, represented by the road traffic upstream internal boundary condition. Finally, the speeds of road traffic and tramway traffic are compared as a criterion to determine whether the tram is blocked by the surrounding road traffic under congestion. The processes of determining internal boundary conditions for road traffic and tramway traffic are demonstrated in Appendix \ref{internal_algorithm_appendix}. For details about determining internal boundary conditions, readers can refer to \cite{xie2026continuum}.

The multimodal traffic internal boundary conditions are incorporated with the upstream and downstream boundary conditions proposed in Section \ref{independent link model} to obtain the cumulative sending and receiving flows, as shown in Equation \ref{multimodal cumulative sending flow} and \ref{multimodal cumulative receiving flow}. 
\begin{equation}
    \overrightarrow{N}^{[m]}(t)=\left\{
    \begin{array}{cc}
         0 & t-\frac{L}{v^{[m]}}<0 \\
         \min\{N^{[m]}(0,t-\frac{L}{v^{[m]}}),\bar{N}^{[m]}_u(L,t),N^{[m]}(L,t-\frac{L}{v^{[m]}})+s^{[m]}\frac{L}{v^{[m]}}\}& t-\frac{L}{v^{[m]}}\geq 0
    \end{array}\right.
    \label{multimodal cumulative sending flow}
\end{equation}
\begin{equation}
    \overleftarrow{N}^{[m]}(t)=\left\{
    \begin{array}{cc}
         \rho_j^{[m]}L & t-\frac{L}{w^{[m]}}<0 \\
         \min\{N^{[m]}(L,t-\frac{L}{w^{[m]}})+\rho_j^{[m]}L,\bar{N}_d^{[m]}(0,t),N^{[m]}(L,t-\frac{L}{w^{[m]}})+s^{[m]}\frac{L}{w^{[m]}}\}& t-\frac{L}{w^{[m]}}\geq 0
    \end{array}\right.
    \label{multimodal cumulative receiving flow}
\end{equation}
\noindent
where $\bar{N}^{[m]}_u(L,t)$ and $\bar{N}^{[m]}_d(0,t)$ are the upstream and downstream internal boundary conditions of traffic modes $m$ at time $t$.

\section{Traffic propagation at nodes}
\label{node model}
This section illustrates the node model of the M-LTM, which achieves traffic transfer between links according to demand constraints, supply constraints, and internal constraints. We discuss two basic types of nodes in urban multimodal traffic networks: tram stops and signalized intersections, given their additional specifications. In the following text, the incoming link of the node is denoted as $i$ or $i'$, while the downstream link is $j$ or $j'$.

\subsection{Tram stop}
\label{tram stop}
The tram dwelling process and dwelling-induced blockage of road traffic are simulated in the tram stop node model. Tramway cumulative flows are determined based on the link demand and supply, and the tram dwelling time, satisfying the First In First Out (FIFO) rule. The road traffic cumulative flows are updated according to whether there are trams at the stop. 

We model the tram stop as a virtual link, whose travel time corresponds to the tram dwelling time and whose capacity is determined by the stop capacity. In this study, the tram dwelling time is predefined according to passenger demand and tram timetables before the tram arrives at the stop. Stochastic terms may be incorporated in future work to capture fluctuations in passenger demand and operational disturbances. In addition, the dwelling time can also be treated as a control variable in tram operation and control problems. The dwelling time of tram arriving at the stop at time $t$ is denoted by $t_d^{[t]}(t)$, and the dwelling time of tram departing from the stop at time $t$ is denoted by $\bar{t}_d^{[t]}(t)$ correspondingly. For example, if the tram arrives at the stop at time $t_1$ is scheduled to dwell for a duration of $\Delta t$, it will depart at time $t_1+\Delta t$. Then, it derives:
\begin{equation}
    t_d^{[t]}(t_1)=\Delta t
\end{equation}
\begin{equation}
    \bar{t}_d^{[t]}(t_1+\Delta t)=\Delta t
\end{equation}
Similar to the link model, the cumulative tramway traffic sending and receiving flows of the tram stop are derived from its cumulative inflow and outflow, as formulated in Equation \ref{tramway stop sending flow} and Equation \ref{tramway stop receiving flow}. Equation \ref{tramway stop sending flow} represents the tram dwelling process, while Equation \ref{tramway stop receiving flow} works when multiple trams dwell at the stop simultaneously.
\begin{equation}
    \overrightarrow{N}^{[t]}_n(t)=\left\{
    \begin{array}{cc}
    0 & t-\bar{t}_d^{[t]}(t)<0\\
    \min \{N_i^{[t]}(L_i,t-\bar{t}_d^{[t]}(t)),N_j^{[t]}(0,t-t_g^{[t]})+1\} & t-\bar{t}_d^{[t]}(t)\geq 0
    \end{array}\right.
    \label{tramway stop sending flow}
\end{equation}
\begin{equation}
    \overleftarrow{N}^{[t]}_n(t)=\left\{
    \begin{array}{cc}
        C^{[t]}_n & t-t_g^{[t]}<0 \\
        \min\{N_j^{[t]}(0,t-t_g^{[t]})+C_n^{[t]},N_i^{[t]}(L_i,t-t_g^{[t]})+1\} & t-t_g^{[t]}\geq 0
    \end{array} \right.
    \label{tramway stop receiving flow}
\end{equation}
\noindent
where $\overrightarrow{N}_n^{[t]}(t)$ and $\overleftarrow{N}_n^{[t]}(t)$ are the cumulative tramway traffic sending and receiving flows of the tram stop at time $t$, and $\bar{t}_d^{[t]}(t)$ is the dwelling time of the tram departing from the stop at time $t$.

The tram stop node model utilizes $\bar{t}_d^{[t]}(t)$ as input, and we provide an approach to obtain this function based on $t_d^{[t]}(t)$ as Equation \ref{tram dwelling time transfer}, while satisfying the FIFO rule. Initially, the function $\bar{t}_d^{[t]}(t)$ is set to the minimum tram headway $t_g^{[t]}$ over the entire simulation horizon, ensuring a minimum departure spacing between consecutive trams at the stop. When a tram arrives at the stop at time $t$ and dwells for $t_d^{[t]}(t)$, the corresponding dwelling duration is assigned to its departure time interval $(t, t+t_d^{[t]}(t))$. Therefore, for each future time $t'\in(t,t+t_d^{[t]}(t)]$, the dwelling time is updated.
\begin{equation}
    \bar{t}_d^{[t]}(t')=\max\{\bar{t}_d^{[t]}(t'),t'-t\},t<t'\leq t+t_d^{[t]}(t)
    \label{tram dwelling time transfer}
\end{equation}
Here, the term $t'-t$ represents the dwelling duration accumulated since the tram arrival at time $t$. The max operator guarantees that the updated dwelling time is not smaller than the existing value, thereby preserving FIFO conditions and preventing later trams from departing earlier than preceding trams when the stop capacity is limited.

\begin{figure}[!h]
    \centering
    \includegraphics[width=0.7\linewidth]{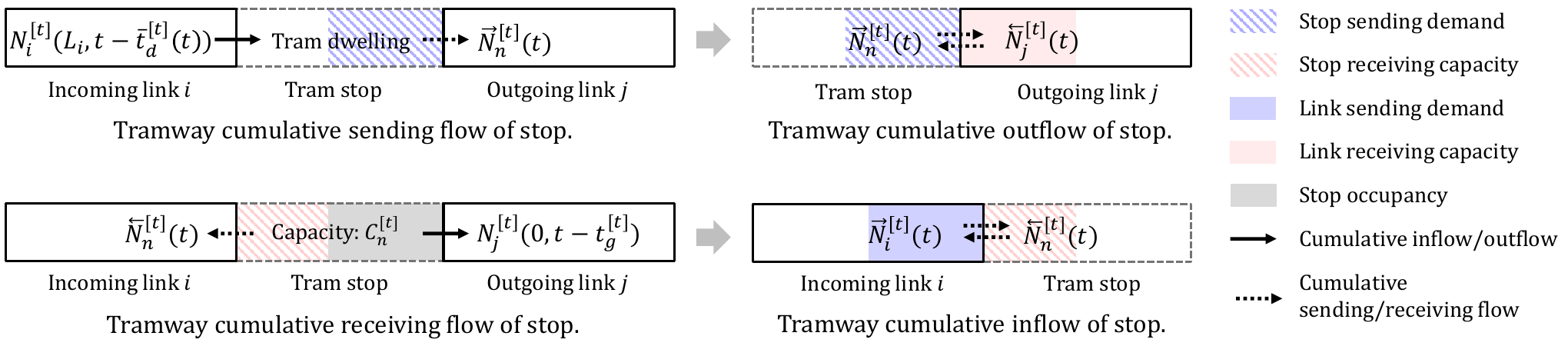}
    \caption{Illustration of tramway traffic transfer at the tram stop node.}
    \label{stop_illustration_tram}
\end{figure}

The cumulative tramway traffic outflow of the incoming link and inflow of the outgoing link are determined by the link and stop constraints, as Equation \ref{tramway traffic cumulative flow of stop incoming link} and Equation \ref{tramway cumulative flow of stop outgoing link}. Figure \ref{stop_illustration_tram} indicates the tramway traffic transfer process considering tram dwellings, the minimum time gap and the stop capacity.
\begin{equation}
    N^{[t]}_i(L_i,t)=\min\{\overrightarrow{N}_i^{[t]}(t),\overleftarrow{N}_n^{[t]}(t)\}
    \label{tramway traffic cumulative flow of stop incoming link}
\end{equation}
\begin{equation}
    N_j^{[t]}(0,t)=\min\{\overleftarrow{N}_j^{[t]}(t),\overrightarrow{N}_n^{[t]}(t)\}
    \label{tramway cumulative flow of stop outgoing link}
\end{equation}

The road traffic dynamics at the tram stop node model simulate the dwelling-induced blockage of road traffic. We consider the unprotected tram stops in the node model. For the protected tram stops where tram dwellings have no impact on the road traffic, the tram stop is modeled as the inhomogeneous node \citep{yperman2005link}. The operation time of road traffic at the stop is the time road traffic has been blocked, which is updated by Equation \ref{waiting time} based on the tram dwelling time. Its initial value is set as $0$ at the beginning, assuming there is no tram at the stop.
\begin{equation}
    \bar{t}_d^{[r]}(t')=\max\{\bar{t}_d^{[r]}(t'),t'-t\} t<t'\leq t+t_d^{[t]}(t)
    \label{waiting time}
\end{equation}
\noindent
where $\bar{t}_d^{[r]}(t)$ is the road traffic waiting time at time $t$.

\begin{figure}[!h]
    \centering
    \includegraphics[width=0.5\linewidth]{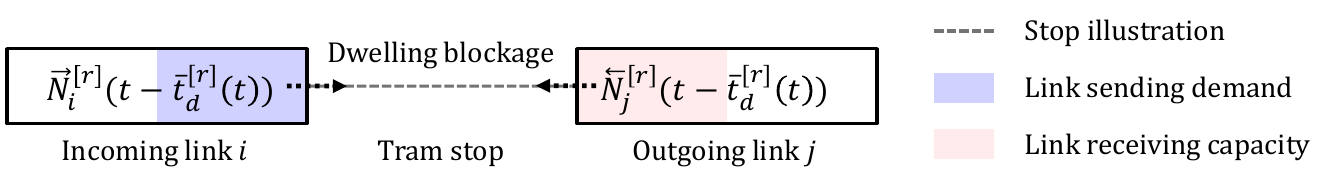}
    \caption{Illustration of road traffic transfer at the tram stop node.}
    \label{stop_illustration_road}
\end{figure}

The road traffic transfer process at tram stops is demonstrated in Figure \ref{stop_illustration_road}. Since the road traffic can not be stored at the stop, the cumulative outflow of the incoming link and inflow of the outgoing link are equal and restricted by the cumulative sending and receiving flows, and the tram stop node constraint as Equation \ref{stop road cumulative flow}.
\begin{equation}
    N_i^{[r]}(L_i,t)=N_j^{[r]}(0,t)=\left\{
    \begin{array}{cc}
    0 & t-\bar{t}_d^{[r]}(t)<0\\
    \min\{\overrightarrow{N}_i^{[r]}(t-\bar{t}_d^{[r]}(t)),\overleftarrow{N}_j^{[r]}(t-\bar{t}_d^{[r]}(t)),N_i^{[r]}(L_i,t-\bar{t}_d^{[r]}(t))+s^{[r]}\bar{t}_d^{[r]}(t)\} & t-\bar{t}_d^{[r]}(t)\geq 0
    \end{array}\right.
    \label{stop road cumulative flow}
\end{equation}

\subsection{General signalized intersection}
\label{intersection without spillback}
The signalized intersection node model simulates the traffic distribution at the intersection regarding the traffic demands from incoming links, supplies of outgoing links, and the signal control. For operational safety, signal control is required at intersections involving trams, and conflicting and merging flows with tramway traffic are prohibited within the same phase. Additionally, tramway and road traffic can be discharged simultaneously within their respective or shared phases without physical interference. In this study, the control scheme is assumed to be given a priori, and its effect is represented by green ratios averaged over cycles \citep{han2015continuum}. In the general signalized intersection, road traffic is prohibited from entering the intersection once its outgoing link reaches capacity. Therefore, road and tramway traffic have no spatio-temporal resource conflicts within the intersection, while the tramway priority is considered in the model. The special cases where road traffic spillback occurs from the outgoing links are addressed separately in Section \ref{signalized intersection with spillback}.

Given the spatial-temporal decoupling of the two traffic modes in the general signalized intersections, we first update the road traffic cumulative flow through an iterative process, and then determine the tramway cumulative flow by a linear program. The road traffic demands are distributed to outgoing links considering the supply constraints. Although the order of distributing outgoing link supplies has no impact on the invariance of the model, it influences the solutions significantly \citep{jabari2016node}. The first considered outgoing links possess more priority. Therefore, we first consider the outgoing links with trams, denoted by set $O^{[t]}$, reflecting the tram priority. To explicitly simulate the complex transient transfer of road traffic, the flow rate is utilized, deriving from the cumulative flow as Equation \ref{raod traffic demand flow} and Equation \ref{road traffic supply flow}.
\begin{equation}
    \overrightarrow{q}_i^{[r]}(t)=\frac{\partial\overrightarrow{N}_i^{[r]}(t)}{\partial t}
    \label{raod traffic demand flow}
\end{equation}
\begin{equation}
    \overleftarrow{q}_j^{[r]}(t)=\frac{\partial\overleftarrow{N}_j^{[r]}(t)}{\partial t}
    \label{road traffic supply flow}
\end{equation}
\noindent
where $ \overrightarrow{q}_i^{[r]}(t)$ is the road traffic demand flow rate of incoming link $i$ at time $t$, and $\overleftarrow{q}_j^{[r]}(t)$ is the road traffic supply flow rate from the outgoing link $j$ at time $t$.

The iterative process identifies the critical outgoing link $j^*$ with the strictest supply constraint in $O^{[t]}$ first, to which the flows are referred as critical flows. If the road traffic demand to $j^*$ is less than its supply capacity, the traffic flow is constrained by the demand and signal control. Otherwise, critical flows will be correspondingly reduced by a reduction index $a_j^{[r]}(t)$, formulated by Equation \ref{reduction index without spillback}. The conservation of turning fraction leads to a reduction in all flows sharing lanes with these critical flows.  This implies that some supply becomes available for other incoming links. Consequently, the unutilized supply can be reallocated to other flows, necessitating an iterative process to distribute the remaining capacity until all supplies or demands are fully exhausted.
\begin{equation}
    a_j^{[r]}(t)=\min\{1,\frac{\overleftarrow{q}_j^{[r]}(t)}{\sum_{i\in I_j}\min\{e_{ij}^{[r]}\overrightarrow{q}_i^{[r]}(t),s_{ij}^{[r]}\}}\}
    \label{reduction index without spillback}
\end{equation}
\noindent
where $a_j^{[r]}(t)$ is the reduction index of outgoing link $j$ at time $t$, $s_{ij}^{[r]}$ is the saturation flow for road traffic from incoming link $i$ to outgoing link $j$, and $I_j$ is the incoming link set of the outgoing link $j$. The saturation flow $s_{ij}^{[r]}$ is determined by the layout of the incoming link $i$ and the green ratio. It will be a changing variable in control optimization problems when the green ratios are control variables. The turning ratio is denoted as $e^{[r]}_{ij}$, formulated as time-independent here. Nevertheless, it can be extended to a time-varying parameter according to the real-time traffic data.

Tramway traffic transfer is updated simultaneously with road traffic through a linear program. Its turning information follows the route information, timetable at the upstream stop, and free flow travel time between the stop and the intersection, which is demonstrated by a binary turning indicator $e^{[t]}_{ij}(t)$. The indicator is a priori, which turns to 1 once the tram is scheduled to pass the intersection. It can be a control variable for signal control and tram scheduling problems. Delays resulting from exogenous disturbances are treated as input in the updated timetable, while the delays caused by multimodal traffic interactions on links are naturally absorbed by the model. Because the tram operations follow the FIFO rule, Equation \ref{intersection conservation law} ensures that the overall tramway tuning ratios from each incoming link are simply shifted backwards in time, preserving the flow conservation at intersections. The linear program consists of the demand and supply constraints and the conservation law.
\begin{equation}
    \max \{ N_i^{[t]}(L_i,t), N_j^{[t]}(0,t) \}
    \label{intersection objective}
\end{equation}
s.t.
\begin{equation}
    N_i^{[t]}(L_i,t)\leq \overrightarrow{N}_i^{[t]}(t), i\in I
    \label{intersection demand constraint}
\end{equation}
\begin{equation}
    N_j^{[t]}(0,t)\leq \overleftarrow{N}_j^{[t]}(t), j\in O
    \label{intersection supply constraint}
\end{equation}
\begin{equation}
    N_j^{[t]}(0,t)=\sum_{i\in\{i'|\int_t{e_{i'j}^{[t]}(t)>0\}}}{\frac{\int_te_{ij}^{[t]}(t)}{\int_t\sum_{j'}e_{ij'}^{[t]}(t)}N_i^{[t]}(L_i,t)}, j\in O
    \label{intersection conservation law}
\end{equation}
Equation \ref{intersection demand constraint} and Equation \ref{intersection supply constraint} regard the demand constraint of incoming links and the supply of the outgoing links, respectively. The conservation law is achieved by Equation \ref{intersection conservation law}.

The process of updating traffic flow at intersections is summarized by Algorithm \ref{flow update without spillback}. The solution of this algorithm satisfies the invariance and holding-free principles defined in \cite{jabari2016node}, which is proved in Appendix \ref{IHF proof}.

\begin{algorithm}[H]
\small
    \label{flow update without spillback}
    \caption{Algorithm of traffic flow update at general signalized intersection.}
    \KwIn{Road traffic cumulative sending flows of incoming links $\overrightarrow{N}^{[r]}_i(t),i\in I$, road traffic cumulative receiving flow of outgoing links $\overleftarrow{N}_j^{[r]}(t),j\in O$}
     \KwIn{Tramway traffic cumulative sending flows of incoming links $\overrightarrow{N}^{[t]}_i(t),i\in I$, tramway traffic cumulative receiving flow of outgoing links $\overleftarrow{N}_j^{[t]}(t),j\in O$}
     \KwIn{Road traffic turning ratio $e_{ij}^{[r]},i\in I,j\in O$, tramway traffic turning ratio $e_{ij}^{[t]}(t),i\in I,j\in O$, road traffic saturation flow $s_{ij}^{[r]},i\in I,j\in O$}
    \KwOut{Road traffic cumulative outflows of incoming links $N^{[r]}_i(L_i,t),i\in I$ and cumulative inflow of outgoing links $N_j^{[r]}(0,t),j\in O$, tramway traffic cumulative outflows of incoming links $N^{[t]}_i(L_i,t),i\in I$ and cumulative inflow of outgoing links $N_j^{[t]}(0,t),j\in O$}
    \textbf{Step 1: Initialize.}\\
    Calculate the road traffic demand and supply flow rate by Equation \ref{raod traffic demand flow} and Equation \ref{road traffic supply flow}.\\
    \textbf{Step 2: Identify the critical outgoing link with tramway traffic.}\\
    Calculate the reduction index $a_j^{[r]}(t)$ for $j\in O^{[t]}$ by Equation \ref{reduction index without spillback}.\\
    $a_{j*}=\min_{j\in O^{[t]}}a_j^{[r]}(t)$\\
    $j^*=\arg\min_{j\in O^{[t]}}a_j^{[r]}(t)$\\
    \If{$a_{j*}\geq 1$}{
    \For{$i\in I$}{
    $q_{ij}^{[r]}=\min\{e_{ij}\overrightarrow{q}_i^{[r]}(t),s_{ij}^{[r]}\}$
    }
    Turn to Step 6.}
    \textbf{Step 3: Reduce the road traffic flow of incoming links to $j^*$.}\\
   Distribute the traffic demands according to the fraction of flow from $i$ to the total demand to the outgoing link $j$:
    $\alpha_{ij}=\frac{e^{[r]}_{ij}}{\sum_{\{i'|j\in O_{i'}\}}e_{i'j}^{[r]}}$\\
    $\bar{q}_{ij}^{[r]}(t)=\min\{\alpha_{ij}\overrightarrow{q}^{[r]}_i(t),s_{ij}^{[r]}\},i\in I,j\in O$\\
    Update the road traffic flow from incoming links of $j^*$.\\
    \For{$i\in I_{j*}$}{
    $q_i^{[r]}(t)=\min\{\min_{j\in O_i}\frac{\overleftarrow{q}_j^{[r]}(t)}{\alpha_{ij}},\sum_{j\in I_j}\bar{q}_{ij}^{[r]}(t)\}$\\
    $q_{ij}^{[r]}(t)=e_{ij}q_i^{[r]}(t)$
    }
    \textbf{Step 4: Update the supply flow of outgoing links.}\\
    \For{$i \in I_{j*}$}{
    \For{$j\in O_i$}{
    $\overleftarrow{q}^{[r]}_j(t)\leftarrow \overleftarrow{q}^{[r]}_j(t)-q_{ij}^{[r]}(t)$\\
    $e_{ij}^{[r]}\leftarrow 0$
    }
    }
    \textbf{Step 5: Return to Step 2.}\\
    \textbf{Step 6: Identify the critical outgoing link without tramway traffic.}\\
    Calculate the reduction index $a_j^{[r]}(t)$ for $j\in O\setminus O^{[t]}$ by Equation \ref{reduction index without spillback}.\\
    $a_{j*}=\min_{j\in O\setminus O^{[t]}}a_j^{[r]}(t)$\\
    $j^*=\arg\min_{j\in O\setminus O^{[t]}}a_j^{[r]}(t)$\\
    \If{$a_{j*}\geq 1$}{Turn to Step 8.}
    \textbf{Step 7: Distribute the remaining road traffic flow.}\\
    Process Step 3 and Step 4.\\
    Return to Step 6.\\
    \textbf{Step 8: Update road traffic cumulative flows.}\\
     $N_i^{[r]}(L_i,t)=\sum_j\int_0^tq_{ij}^{[r]}(\tau){\rm d}\tau$\\
    $N_j^{[r]}(0,t)=\sum_i\int_0^tq_{ij}^{[r]}(\tau){\rm d}\tau$\\
    \textbf{Step 9: Update tramway traffic cumulative flows.}\\
    Solve the program of Equation \ref{intersection objective} with Equation \ref{intersection demand constraint}, Equation \ref{intersection supply constraint}, and Equation \ref{intersection conservation law}.\\    
    
\end{algorithm}

\subsection{Signalized intersection with spillback}
\label{signalized intersection with spillback}
Road traffic is prohibited from leaving incoming links once the outgoing link inflow reaches its supply capacity in traditional intersection models, as described in Section \ref{intersection without spillback}. However, when the outgoing link becomes severely congested, road traffic may spill back into the intersection, blocking conflicting road and tramway flows and causing intersection gridlock. To address this challenge, we extend the general signalized intersection node model in the previous section. The physical space of the intersection is conceptualized as a buffer supply extension of the downstream outgoing link. This enables the model to simulate the temporary road traffic queuing within the intersection. Building upon this, the spillback is identified, and the blockage of conflicting road and tramway traffic is modeled.

First, the original reduction index is modified by incorporating the intersection buffer supply capacity, which holds the spilling back road traffic, as Equation \ref{reduction index with spillback}.
\begin{equation}
    a_j^{[r]}(t)=\min\{1,\frac{\overleftarrow{q}_j^{[r]}(t)+\tilde{C}^{[r]}}{\sum_{i\in I_j}\min\{e_{ij}^{[r]}\overrightarrow{q}_i^{[r]}(t),s_{ij}^{[r]}\}}\}
    \label{reduction index with spillback}
\end{equation}
\noindent
where $\tilde{C}^{[r]}$ is the intersection buffer supply capacity of road traffic. It measures the additional flow rate that the intersection can temporarily hold. Under spillback conditions, the total supply capacity for road traffic to outgoing link $j$ is no longer bounded solely by the link supply $\overleftarrow{q}_j^{[r]}(t)$, but is augmented by the intersection buffer supply capacity $\tilde{C}^{[r]}$. Therefore, the aggregate road traffic demand is compared with this extended value to determine whether the supply satisfies the demand.

Them, the binary indicator $\tilde{\varepsilon}_j(t)$ is introduced to suggest whether the intersection suffers a spillback from outgoing link $j$, which implies the blockage of other conflicting road and tramway traffic before updating their cumulative flow.
\begin{equation}
   \tilde{\varepsilon}_j(t)=\left\{
    \begin{array}{cc}
    1 & \text{the outgoing link }j\text{ spills back.}\\
    0 & \text{otherwise.}
    \end{array}\right.
    \label{varepsilon}
\end{equation}
Once the spillback occurs $\tilde{\varepsilon}_j(t)=1$, other conflicting traffic is prohibited from entering the intersection until the congestion dissipates. The road traffic demand to the congested outgoing link $j$ is the queuing traffic at the intersection, as shown in Equation \ref{spillback demand}. The spillback disappears when all the remaining traffic enters link $j$.
\begin{equation}
    \sum_{i\in I_j}\overrightarrow{q}^{[r]}_{ij}(t)=\frac{\partial(\sum_{i\in I_j}e_{ij}^{[r]}N^{[r]}_i(L_i,t)-N^{[r]}_j(0,t))}{\partial t}
    \label{spillback demand}
\end{equation}
\begin{figure}[!h]
    \centering
    \includegraphics[width=\linewidth]{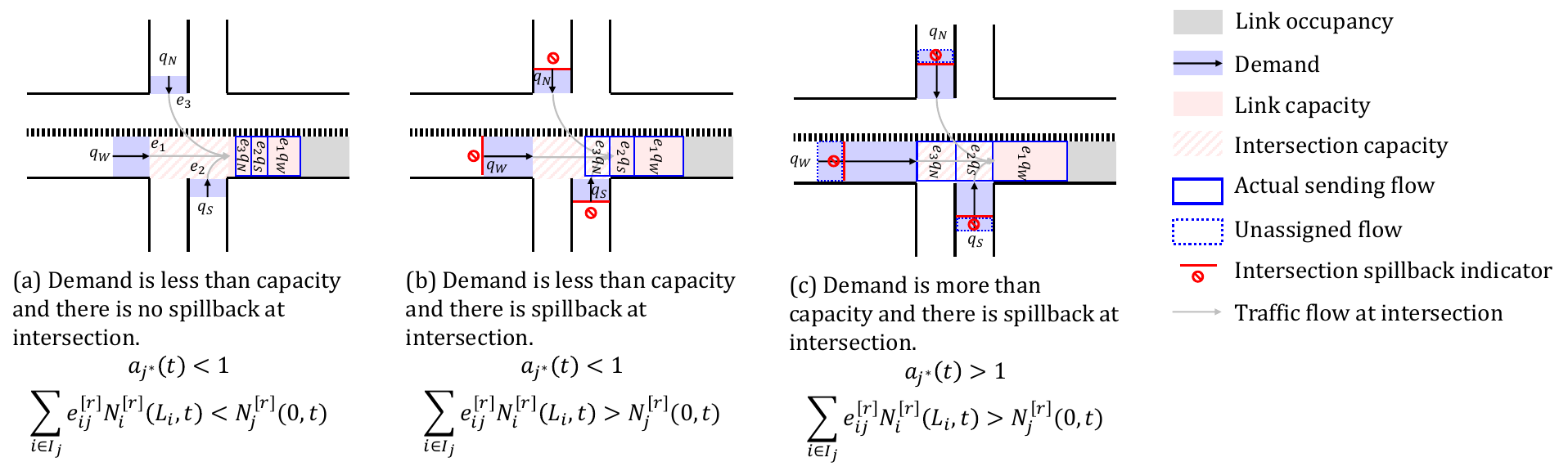}
    \caption{Illustration of road traffic distribution in the signalized intersection node model with spillback, where $a_j$ is defined in Equation \ref{reduction index with spillback}. The turning ratios $e_1$,  $e_2$, and $e_3$ are omitted in Figure (b) and (c).}
    \label{spill_back_intersection}
\end{figure}
Figure \ref{spill_back_intersection} illustrates three possible situations at the signalized intersection with the consideration of spillback inside. When the road traffic demand is less than the outgoing link capacity, the traffic distribution is identical to that of the general node model, as shown in Figure \ref{spill_back_intersection} (a). Once the total demand exceeds the outgoing link capacity, the road traffic spills back into the intersection and causes intersection gridlock. If the intersection buffer supply capacity can hold the exceeding demand, all road traffic demand directed toward this outgoing link is satisfied while other conflicting traffic flows are blocked, as Figure \ref{spill_back_intersection} (b). Figure \ref{spill_back_intersection} (c) indicates that once the traffic demand is more than the extended capacity with the intersection buffer supply, the traffic demand to this outgoing link is reduced in proportion as well. The feasible flow discharged from each incoming link to this critical outgoing link is scaled down according to its turning ratio fraction, similar to the algorithm of the general intersection node mode. This implies that an incoming link contributing a higher share of the turning ratio toward the critical outgoing link will experience a proportionally larger absolute reduction, thereby maintaining the invariance principle during capacity allocation.

The road traffic spillback also blocks tramway traffic and leads to a delay at the intersection. This is modeled by the modified intersection conservation law as Equation \ref{modified intersection conservation law}. If there is no spillback, the cumulative flows are distributed as scheduled, while the tram can not pass the intersection as the timetable under the spillback condition. The cumulative tramway traffic flows are updated by the program of Equation \ref{intersection objective} with Equation \ref{intersection demand constraint}, Equation \ref{intersection supply constraint}, and Equation \ref{modified intersection conservation law}.
\begin{equation}
    N_j^{[t]}(0,t)=\sum_{i,\int_t{e_{ij}^{[t]}(t)>0}}{\frac{\int_te_{ij}^{[t]}(t)-\tilde{\varepsilon}_j(t)e_{ij}^{[t]}}{\int_t\sum_{j'}e_{ij'}^{[t]}(t)-\sum_{j'}\tilde{\varepsilon}_{j'}(t)e^{[t]}_{ij}}N_i^{[t]}(L_i,t)}, j\in O
    \label{modified intersection conservation law}
\end{equation}

The process of updating traffic cumulative flow at signalized intersections concerning spillback is shown in Algorithm \ref{flow update with spillback}, in which the spillback is identified first. The invariance and holding-free principles are satisfied since the algorithm in this section remains the general framework of Algorithm \ref{flow update without spillback}.\\
\begin{algorithm}[p]
\small
    \caption{Algorithm of traffic flow update at signalized intersection considering spillback.}
    \label{flow update with spillback}
    \KwIn{Road traffic cumulative sending and receiving flows$\overrightarrow{N}^{[r]}_i(t),i\in I$, $\overleftarrow{N}_j^{[r]}(t),j\in O$, spillback indicator $\tilde{\varepsilon}_j(t)$}
    \KwIn{Tramway traffic cumulative sending and receiving flows$\overrightarrow{N}^{[t]}_i(t),i\in I$, $\overleftarrow{N}_j^{[t]}(t),j\in O$}
     \KwIn{Road traffic turning ratio $e_{ij}^{[r]},i\in I,j\in O$, tramway traffic turning ratio $e_{ij}^{[t]}(t),i\in I,j\in O$, road traffic saturation flow $s_{ij}^{[r]},i\in I,j\in O$, intersection buffer supply capacity $\tilde{C}$}
    \KwOut{Road traffic cumulative outflows and inflows $N^{[r]}_i(L_i,t),i\in I$, $N_j^{[r]}(0,t),j\in O$, tramway traffic cumulative outflows and inflows $N^{[t]}_i(L_i,t),i\in I$, $N_j^{[t]}(0,t),j\in O$}
    \textbf{Step 1: Initialize.}\\
    Calculate the road traffic demand and supply flow by Equation \ref{raod traffic demand flow} and Equation \ref{road traffic supply flow}.\\
    Initialize the candidate incoming and outgoing link sets which are not blocked: $\tilde{I}=I$, $\tilde{O}=O$\\
    \textbf{Step 2: Identify spillback phenomenon.}\\
    \For{$j\in O$}{
    \If{$\tilde{\varepsilon}_j(t)=1$}{
    \For{$i$ influenced by $j$ ($j$ is at the downstream of $i$ or traffic flow from $i$ has conflicts with those to $j$)}{
    $e_{ij'}^{[r]}\leftarrow 0$, 
    $i \setminus \tilde{I}$,
    $q^{[r]}_{ij'}=0\ \forall j'\in O$.\\
     $N_i^{[r]}(L_i,t)=\sum_{j'\in O}\int_0^tq_{ij'}^{[r]}(\tau){\rm d}\tau$
    }
    $N_j^{[r]}(0,t)=\max\{\overleftarrow{N}^{[r]}_j(t),\sum_{i\in I_j}e^{[r]}_{ij}N^{[r]}_i(L_i,t)\}$, $j\setminus \tilde{O}$.\\
    }
    }
    \textbf{Step 3: Identify the critical outgoing link without spillback blockage.}\\
   Calculate the reduction index $a_j^{[r]}(t)$ for $j\in O^{[t]}\cap\tilde{O}$ by Equation \ref{reduction index with spillback}.\\
    $a_{j*}=\min_{j\in O^{[t]}\cap\tilde{O}}a_j^{[r]}(t)$, $j^*=\arg\min_{j\in O^{[t]}\cap\tilde{O}}a_j^{[r]}(t)$\\
    \If{$a_{j*}\geq 1$}{Turn to Step 7.}
    \textbf{Step 4: Reduce the road traffic flow of incoming links to $j^*$.}\\
Distribute the traffic demands according to the fraction of flow from $i$ to the total demand to the outgoing link $j$:
    $\alpha_{ij}=\frac{e^{[r]}_{ij}}{\sum_{\{i'|j\in O_{i'}\}}e_{i'j}^{[r]}}$\\
    $q_{ij}^{[r]}(t)=\min\{\alpha_{ij}\overrightarrow{q}^{[r]}_i(t),s_{ij}^{[r]}\},i\in \tilde{I},j\in \tilde{O}$\\
    \For{$i\in I_{j*}\cap\tilde{I}$}{
    $q_i^{[r]}(t)=\min\{\min_{j\in O_i\cap\tilde{O}}\frac{\overleftarrow{q}_j^{[r]}(t)+\tilde{C}}{\alpha_{ij}},\sum_{j\in I_j}\overrightarrow{q}_{ij}^{[r]}(t)\}$, $q_{ij}^{[r]}(t)=\alpha_{ij}q_i^{[r]}(t)$
    }
    \textbf{Step 5: Update the supply flow of outgoing links.}\\
    \For{$i\in I_{j*}\cap\tilde{I}$}{
    \For{$j\in O_i\cap\tilde{O}$}{
    $\overleftarrow{q}^{[r]}_j(t)\leftarrow \overleftarrow{q}^{[r]}_j(t)-q_{ij}^{[r]}(t)$, $e_{ij}^{[r]}\leftarrow 0$
    }
    }
    \textbf{Step 6: Return to Step 3.}\\
    \textbf{Step 7: Identify the critical outgoing link without tramway traffic.}\\
    Calculate the reduction index $a_j^{[r]}(t)$ for $j\in \tilde{O}\setminus ( O^{[t]}\cap\tilde{O})$ by Equation \ref{reduction index without spillback}.\\
    $a_{j*}=\min_{j\in  \tilde{O}\setminus ( O^{[t]}\cap\tilde{O})}a_j^{[r]}(t)$, $j^*=\arg\min_{j\in \tilde{O}\setminus ( O^{[t]}\cap\tilde{O})}a_j^{[r]}(t)$\\
    \If{$a_{j*}\geq 1$}{Turn to Step 9.}
    \textbf{Step 8: Distribute the remaining road traffic flow.}\\
    Process Step 4 and Step 5.\\
    Return Step 7.\\
    \textbf{Step 9: Update road traffic cumulative flows.}\\
     $N_i^{[r]}(L_i,t)=\sum_j\int_0^tq_{ij}^{[r]}(\tau){\rm d}\tau,i\in \tilde{I}$\\
    $N_j^{[r]}(0,t)=\sum_j\int_0^tq_{ij}^{[r]}(\tau){\rm d}\tau,j\in \tilde{O}$\\
    \textbf{Step 10: Update tramway traffic cumulative flows.}\\
    Solve the program of Equation \ref{intersection objective} with Equation \ref{intersection demand constraint}, Equation \ref{intersection supply constraint}, and Equation \ref{modified intersection conservation law}.\\
    \textbf{Step 11: Update the spillback indicator.}\\
    \For{$j\in O$}{
    \If{$\sum_{i\in I_j}e^{[r]}_{ij}N^{[r]}_i(L_i,t)>N_j^{[r]}(0,t)$}{
    $\tilde{\varepsilon}_j(t')=1,t'>t$
    }
    \Else{
    $\tilde{\varepsilon}_j(t')=0,t'>t$
    }
    }
\end{algorithm}

\section{Case study}
\label{case study}
Three case studies are conducted in this section to demonstrate the ability of the proposed model in multimodal traffic analysis and control. The first case study in Section \ref{one node} illustrates the traffic transfer processes captured by the node model at a tram stop network and a signalized intersection network. Then, the proposed model is implemented on a synthetic artery with three stops and seven intersections to demonstrate the multimodal traffic dynamics and delay propagation in Section \ref{arterial network}. Finally, we apply the proposed model to the Dresden traffic network and compare the results with real-world data to verify the accuracy and practicality of the model in Section \ref{dresden_case}.

The GEH statistic and normalized deviation are introduced to evaluate the accuracy of simulation results of the proposed model compared to baselines. The GEH statistic is widely applied in traffic modeling to compare two sets of traffic volumes, whose value is less than 5 of an acceptable model  \citep{unit2013m3}. The GEH statistic is calculated by Equation \ref{GEH}.
\begin{equation}
    GEH=\sqrt{\frac{2(q^{[r]}_{model}-q^{[r]}_{baseline})^2}{q^{[r]}_{model}+q^{[r]}_{baseline}}}
    \label{GEH}
\end{equation}
\noindent
where $GEH$ is the GEH statistic, and $q_{model}^{[r]}$ and $q_{baseline}^{[r]}$ are the average road traffic flow rate (veh/h) of the simulation duration. The baseline is the SUMO simulation results in Section \ref{one node} and Section \ref{arterial network}, while it is the real-world data in Section \ref{dresden_case}.

While the GEH statistic provides an overall evaluation of the model accuracy during the simulation time, it may smooth out cumulative deviations in traffic dynamics. To further validate the model's ability to reproduce the traffic evolution, the normalized deviation is employed, as shown in Equation \ref{deviation}.

\begin{equation}
    D=\frac{\sqrt{\frac{1}{T}\int_0^T(N^{[r]}_{model}(t)-N^{[r]}_{baseline}(t))^2dt}}{Q}
    \label{deviation}
\end{equation}
\noindent
where $D$ is the normalized deviation, $N^{[r]}_{model}(t)$ and $N^{[r]}_{baseline}(t)$ are the cumulative road traffic flows at time $t$ obtained from the model and baseline, $T$ is the total simulation time, and $Q$ is the total cumulative road traffic flow.

\subsection{One-node network}
\label{one node}
We first illustrate the traffic evolution simulated by the proposed node model at a tram stop and a signalized intersection, respectively. The tram dwelling process, dwelling-induced congestion, and road traffic spillback at intersection are exhibited in the results, and the simulation results are compared with those of the microscopic simulator SUMO (Simulation of Urban MObility).

Parameters of the multimodal traffic in this case study are shown in Table \ref{parameter_1}. The link comprises two lanes with a length of 320 m, with additional dedicated tracks for tramway traffic. The total simulation time is 3600 s. The road traffic demand is 720 veh/h during the first 1200 s, followed by 1440 veh/h during the next 1200 s, and then decreases to 1080 veh/h.
\begin{table}[H]
    \centering
    \begin{tabular}{c|c|c|c|c|c|c}\hline
         $v_f^{[r]}$&  $v_f^{[t]}$&  $w^{[r]}$&  $L_t$&  $t_g$&  $\rho_j^{[r]}$& $s^{[r]}$ \\\hline
         14 $m/s$&  8 $m/s$&  5.6 $m/s$&  30 $m$&  10 $s$&  125 $veh/km/lane$& 1800 $veh/h/lane$\\ \hline
    \end{tabular}
    \caption{Parameters of multimodal traffic flow.}
    \label{parameter_1}
\end{table}

The tram stop network consists of an unprotected stop with its incoming and outgoing links. The capacity of the stop is designed as one or two in the experiments to compare the influence of stop capacity. Two tram lines with departure intervals of eight minutes (480 s) and six minutes (360 s) are scheduled to dwell at the stop for 100 s and 150 s, respectively. A 120 s delay on the entry of the second tram of Line 2 is simulated to depict the traffic dynamics under disturbances.

\begin{figure}[!h]
    \centering
    \includegraphics[width=0.7\linewidth]{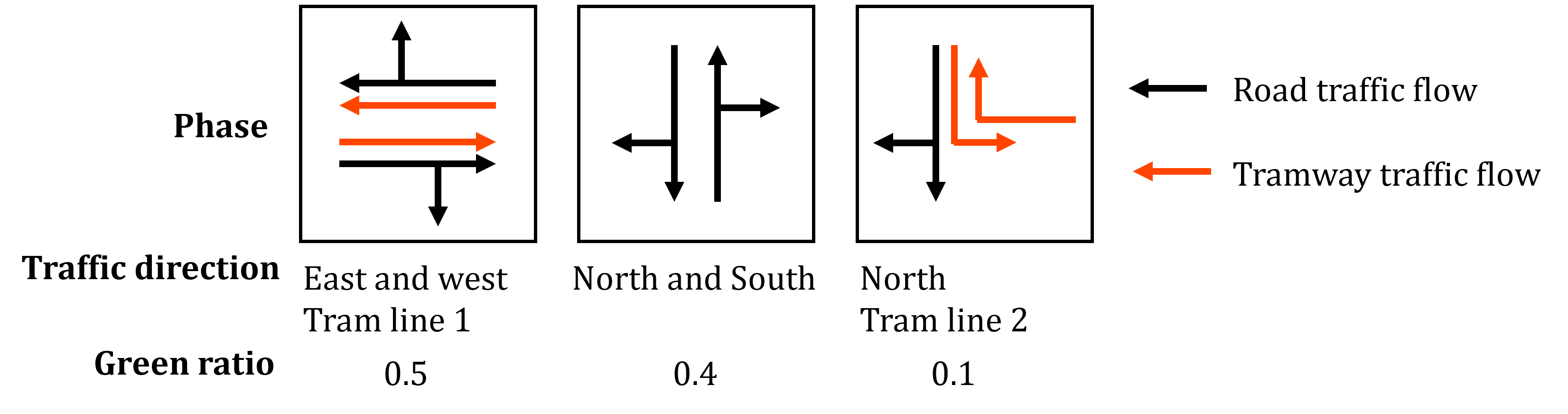}
    \caption{Signal scheme of the intersection network.}
    \label{signal_scheme}
\end{figure}

The intersection network comprises a signalized intersection with four incoming links and four outgoing links. The incoming links are divided into one through lane and one through-right turn lane. The turning ratio of through and right turn flow is 0.7 and 0.3. The signal scheme is depicted in Figure \ref{signal_scheme}. Additionally, two traffic signals with a green ratio of 0.3 are implemented at the end of the south and east outgoing links to simulate bottlenecks, while that at the end of the east outgoing link is removed after the 2400 s. Two bidirectional tram lines are considered. The first tram line passes through the north and east links with a frequency of eight minutes (480 s), while the second one passes through the west and east links with a frequency of six minutes (360 s). The two tram lines merge or diverge at the intersection. More priorities are given to the outgoing links with tramway traffic (north, west, and east outgoing links).

\subsubsection{Tram stop network}
We analyze multimodal traffic dynamics at the tram stop with dedicated ROWs on incoming and outgoing links to demonstrate the tram dwelling process and dwelling-induced congestion on road traffic at the tram stop node model.

\begin{figure}[!h]
    \centering
    \includegraphics[width=\linewidth]{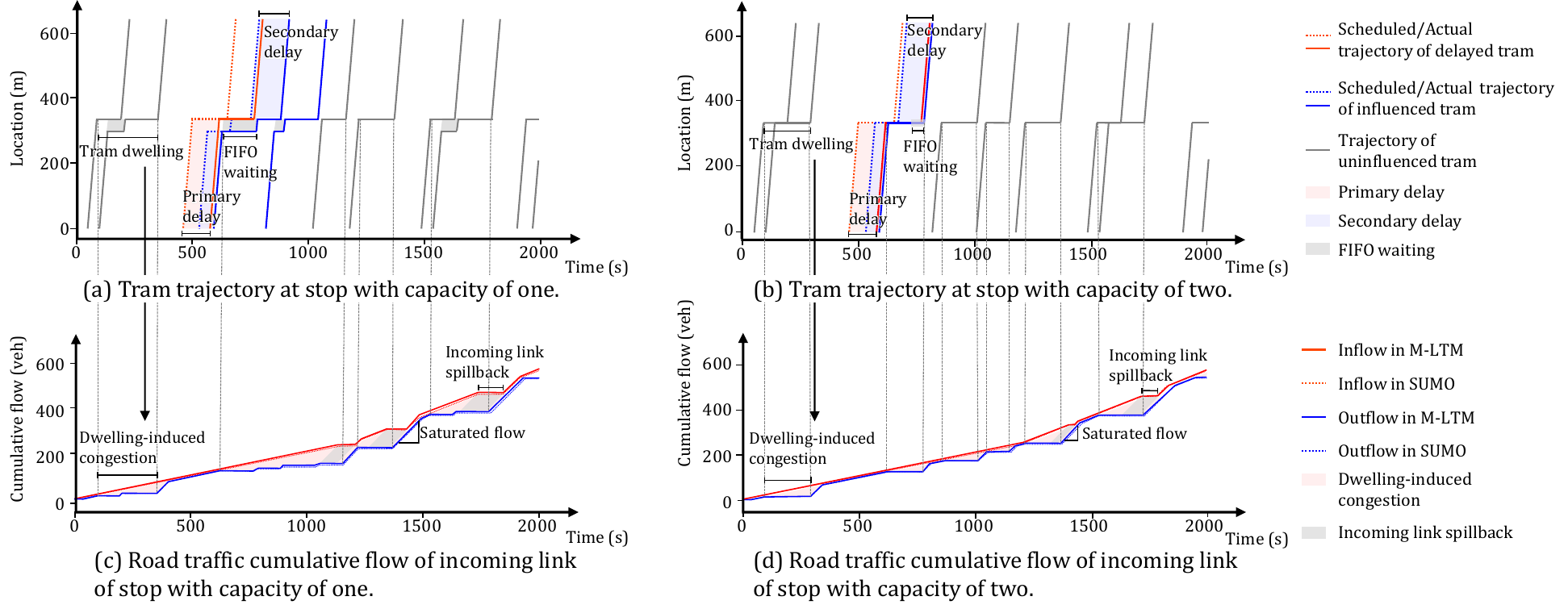}
    \caption{Tram trajectories and road traffic cumulative flows of the incoming link in the tram stop network}
    \label{tram_stop_results}
\end{figure}

The tramway and road traffic dynamics during the first 2000 s are depicted by the tram trajectories and road traffic cumulative flows of the incoming link in Figure \ref{tram_stop_results}. The tram trajectories are derived from the cumulative flows to show the tram operations in detail, assuming the tramway traffic speed is constant on links. Results indicate that a stop capacity of two is significantly more efficient than a capacity of one. The primary delay on the third tram, denoted by the red trajectories in Figure \ref{tram_stop_results}, causes a secondary entry delay of the subsequent tram on the incoming link. In the two-capacity scenario, tram operations recover in a short period after the delay, as two trams can dwell at the stop at the same time. The fourth tram dwells at the stop for an additional 50 seconds after the previous tram according to the FIFO rule. Conversely, with a stop capacity of one, the delay propagates to the next two trams. The fourth tram has to wait to enter the stop until the previous one finishes dwelling, and this delayed process further influenced the fifth tram.

Two types of congestion are identified through the road traffic cumulative inflow and outflow, respectively. The dwelling-induced congestion is observed from the cumulative outflow. Remarkably longer time of congestion is noticed at the stop with a capacity of one, since trams have to dwell at the stop sequentially and the stop is occupied for a longer time. The congestion caused by long-time blockage before the tram stop propagates to the upstream of the incoming link under heavy road traffic demands. This is reflected by a constant cumulative inflow, indicating that the road traffic spills back from the incoming link. Because the blockage time is longer at the stop with less capacity, this type of congestion is more frequent and persistent at the stop with capacity of one. This suggests that the interactions between the two traffic modes may cause severe congestion when the traffic demands are heavy while the infrastructure capacity is not enough.

The cumulative flows derived from the M-LTM are validated against SUMO. The SUMO simulation results are denoted by the dashed lines in Figure \ref{tram_stop_results}, which are highly consistent with the M-LTM prediction. The normalized deviations and GEH statistics between the two simulations are shown in Table \ref{D_n general stop}. The deviations are less than 1\%, and the GEH statistics are below 1, confirming the accuracy of the proposed macroscopic model relative to the microscopic simulation.

\begin{table}
    \centering
    \begin{tabular}{ccccccc}\toprule
         &  \multicolumn{2}{c}{Incoming link inflow}&  \multicolumn{2}{c}{Incoming link outflow}&  \multicolumn{2}{c}{Outgoing link outflow}\\
         &  $D$&  $GEH$&  $D$&  $GEH$&  $D$& $GEH$\\\midrule
         Stop capacity of one&  0.39\%&  0.06&  0.42\%&  0.19&  0.63\%& 0.58\\
         Stop capacity of two&  0.36\%&  0.06&  0.46\%&  0.28&  0.69\%& 0.69\\ \bottomrule
    \end{tabular}
    \caption{Normalized deviation and GEH statistic of road traffic cumulative flow of the general M-LTM in tram stop network.}
    \label{D_n general stop}
\end{table}

\subsubsection{Signalized intersection network}
The traffic dynamics at the signalized intersection captured by the two node models in Section \ref{intersection without spillback} and Section \ref{signalized intersection with spillback} are discussed in this section. We first analyze the traffic evolution and tram priority at the intersection using the general intersection node model, followed by the spillback phenomenon under congested conditions.

\begin{figure}[!h]
    \centering
    \includegraphics[width=\linewidth]{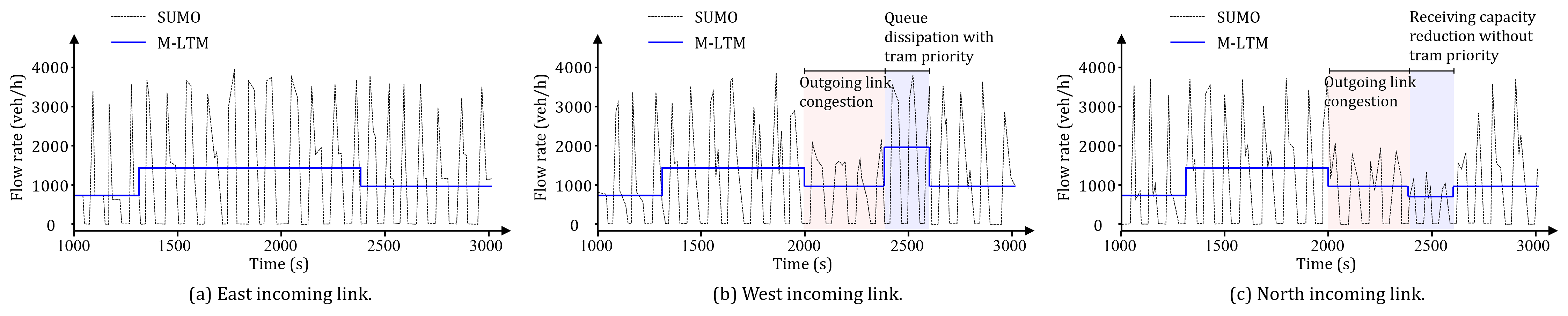}
    \caption{Outflow rate of incoming links in the general intersection node.}
    \label{intersection_result_1}
\end{figure}

\begin{figure}
    \centering
    \includegraphics[width=0.7\linewidth]{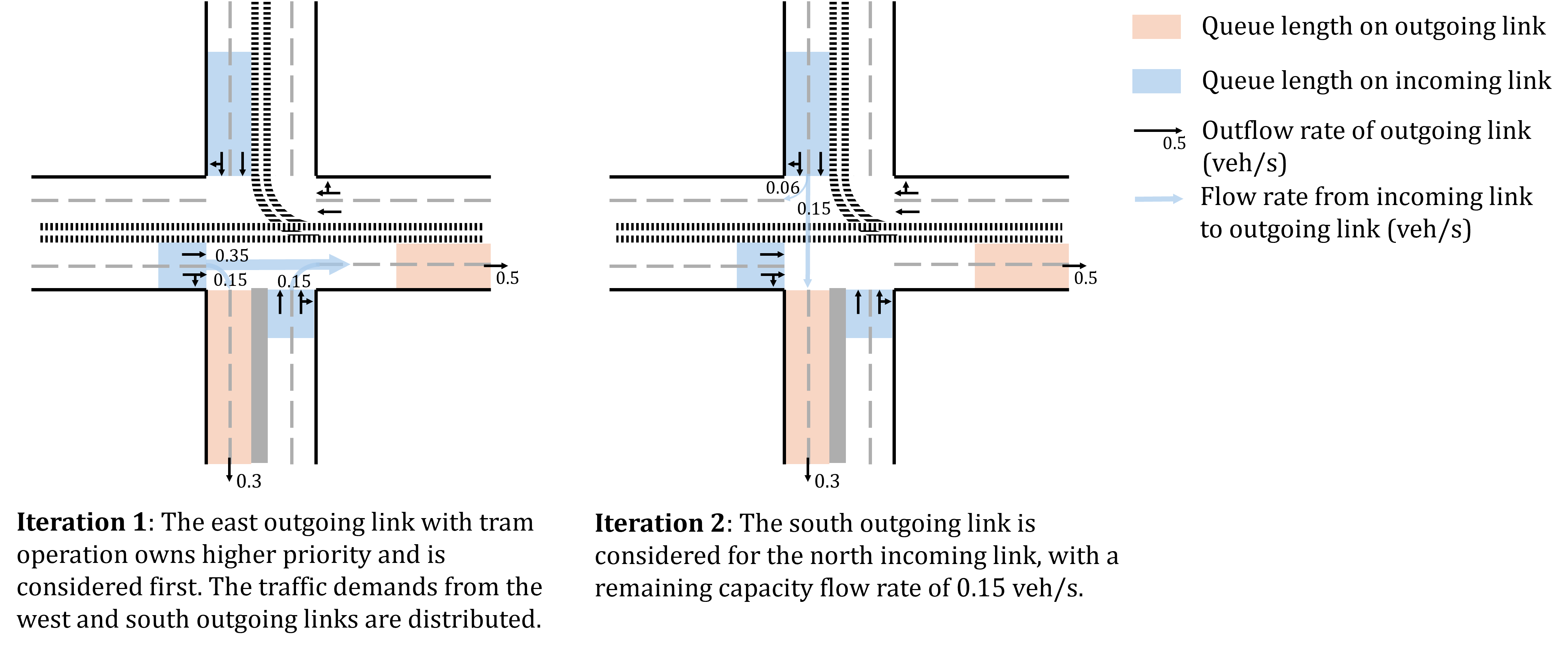}
    \caption{Traffic distribution process at the intersection during 2400 s and 2600 s.}
    \label{intersection_traffice_distribution_illustration}
\end{figure}

The general signal intersection node model without road traffic spillback is implemented first. The road traffic outflow rates of incoming links derived from the cumulative flows during the 1000s and 3000 s are exhibited in Figure \ref{intersection_result_1}. With the unsaturated demand in the first 1200 s, all traffic flows operate without congestion. The bottlenecks at the downstream of the east and south outgoing links lead to congestion when the demand increases. The congestion propagates to the upstream in 2000 s, which restricts the traffic from incoming links entering these two outgoing links and causes a decrease in the flow rates. Once the signal at the end of the east outgoing link is removed, the queue begins to dissipate, and the inflow rate increases to the saturation flow until the queue on incoming links totally dissipates. Although the south outgoing link is more congested with a smaller reduction index $a_j$, a priority is given to the east outgoing link with tram operation. Therefore, the traffic demands from the south and west incoming links are distributed first, followed by the demand from the north incoming link. Since the limited receiving capacity remains after distributing traffic from the west incoming link, the outflow rate of the north incoming link decreases, which aggregates the congestion. The road traffic distribution logic during this period from the 2400s to 2600 s is demonstrated by Figure \ref{intersection_traffice_distribution_illustration}.

The simulation results of the proposed M-LTM and SUMO are quantified and compared by the normalized deviations and GEH statistics, as shown in Table \ref{D_n general intersection}. The main deviations result from the difference between signal representations. The M-LTM utilizes continuous green ratios to approximate the capacity reduction by signal control, while the microscopic simulator SUMO considers the phase time. However, the GEH statistics and normalized deviations remain marginal, verifying that the proposed model can capture the traffic dynamics at the signalized intersection.

\begin{table}
    \centering
    \begin{tabular}{ccccccccc}\toprule
         &  \multicolumn{2}{c}{North}&  \multicolumn{2}{c}{South}&  \multicolumn{2}{c}{West}&  \multicolumn{2}{c}{East}\\
         &  $D$&  $GEH$&  $D$&  $GEH$&  $D$&  $GEH$&  $D$& $GEH$\\\midrule
         Downstream of incoming link&  2.31\%
&  0.89&  2.03\%
&  0.78&  1.89\%
&  0.65&  1.56\%
& 0.66\\
         Upstream of outgoing link&  2.41\%&  0.92&  4.36\%&  1.32&  2.36\%&  0.83&  3.57\%& 1.36\\\bottomrule
    \end{tabular}
    \caption{Normalized deviation and GEH statistic of road traffic cumulative flow of the general LTM in signalized intersection network.}
    \label{D_n general intersection}
\end{table}

\begin{figure}[!h]
    \centering
    \includegraphics[width=0.9\linewidth]{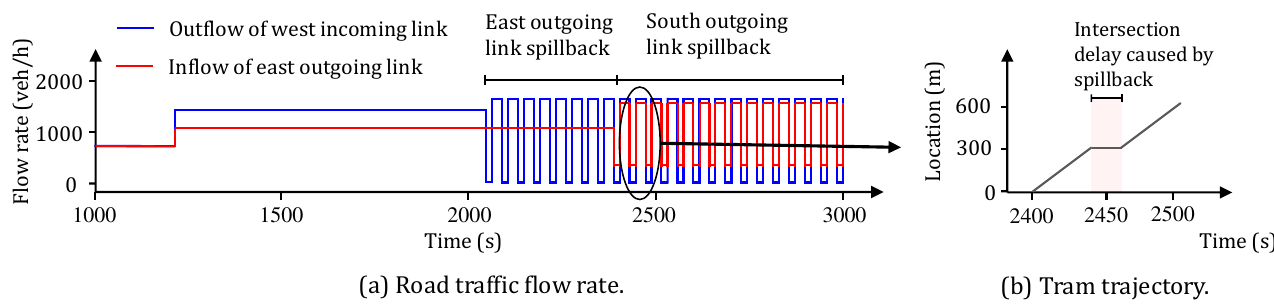}
    \caption{Road traffic flow rate and tram trajectory at the intersection with spillback.}
    \label{intersection_spillback}
\end{figure}

The road traffic spillback under congestion is further simulated using the model in Section \ref{signalized intersection with spillback}. The intersection buffer supply capacity in the model is three vehicles. Figure \ref{intersection_spillback} (a) depicts the road traffic outflow rate of the west incoming link and inflow rate of the east outgoing link during 1000 s and 3000 s. When the east and south outgoing links are congested, vehicles enter the intersection and induce gridlock. Although road traffic from the east incoming link is not influenced by the congestion on outgoing links under general conditions, it is also blocked by the congested vehicles at the intersection in this model. Once the bottleneck at the downstream of the east outgoing link is removed, the receiving capacity of this link recovers. The south link is identified to spill back to the intersection, and the traffic flow to this outgoing link is distributed first. As a result, the outflow of the west incoming link is still blocked during spillback and reduced when the intersection is empty. The outflow oscillation and reduction further lead to the fluctuation of the inflow of the east outgoing link.
The simulation results suggest that the intersection spillback not only reduces the inflow and outflow of links with bottlenecks and their incoming links, but also leads to severe congestion on other links and gridlock at the intersection. Moreover, the frequent stop-and-go at the intersection causes traffic instability, and the congestion may propagate to the whole network.

The spillback of road traffic at the intersection also delays tramway traffic and induces a total increase of 132 s in its operation time. Before the 1200 s, trams operate as scheduled in the two models, which is consistent with the conclusion from the results of road traffic. With the increase in road traffic demand during the 1200 s and 2400 s, congestion appears and propagates to the upstream of the east outgoing link. After the 2000 s, the east outgoing link reaches its capacity. The congestion propagates to the intersection and incoming links, and contributes to extra tramway delays. Trams are delayed to enter the outgoing link, once the road traffic spillback is considered. An example of tram trajectory during this period is provided in Figure \ref{intersection_spillback} (b). After the bottleneck at the downstream of the outgoing link is removed, the tramway traffic from the west incoming link recovers in a short time, since the congestion dissipates within 200 s. Nevertheless, trams from the north incoming link are still delayed, which resulted from the limited receiving capacity of the south outgoing link. Given the spillback from the south outgoing link, the tramway traffic is blocked by the queue from the south outgoing link, which causes additional delays at the intersection. The results demonstrate that road traffic spillback at the intersection has significant impacts not only on road traffic but also on tramway traffic, as the congestion blocks trams from entering the intersection, even with their dedicated ROW.

\subsection{Arterial network}
\label{arterial network}
The proposed multimodal traffic dynamic model is implemented on a synthetic arterial network in this section to illustrate its ability to simulate traffic evolution with multimodal interactions and tram delay propagation.

\begin{figure}
    \centering
    \includegraphics[width=0.9\linewidth]{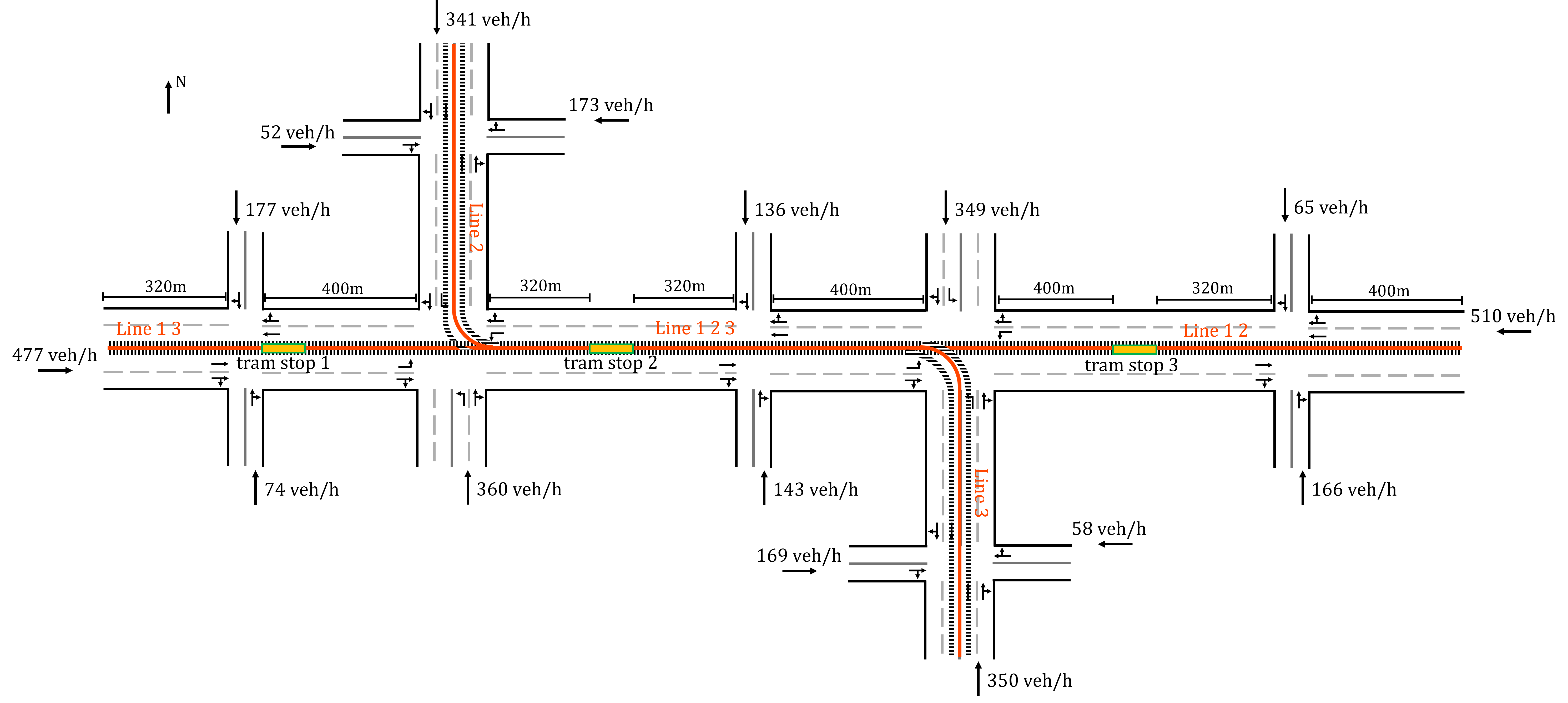}
    \caption{Layout of the arterial network.}
    \label{arterial layout}
\end{figure}

The network consists of seven signalized intersections and three tram stops with capacities of one, as shown in Figure \ref{arterial layout}. The horizontal artery comprises two lanes for road traffic and a dedicated tram track in each direction, while the two traffic modes share the ROW on longitudinal arteries with inter-modal interactions. The fundamental parameters are the same as those in the one-node networks (as shown in Table \ref{parameter_1}), and the road traffic relative flow with moving bottleneck is 1350 veh/h. The road traffic demand of the whole network is 3600 veh/h, and the exact demands at source nodes are shown in Figure \ref{arterial layout}. Three bidirectional tram routes are considered in the network, merging or diverging on the horizontal artery. The departure interval of Line 1 is five minutes (300 s), while the departure intervals for the other tram lines are ten minutes (600 s). Tram stop 1 is next to the intersection, which is regarded as a protected stop and modeled together with the intersection in M-LTM. The other two stops are unprotected stops, in which the dwelling-induced congestion of road traffic is captured. The dwelling time at Stop 1 is 60 s, and that of the other stops is 90 s.

\begin{figure}[!h]
    \centering
    \begin{subfigure}[b]{0.49\linewidth}
        \includegraphics[width=\linewidth]{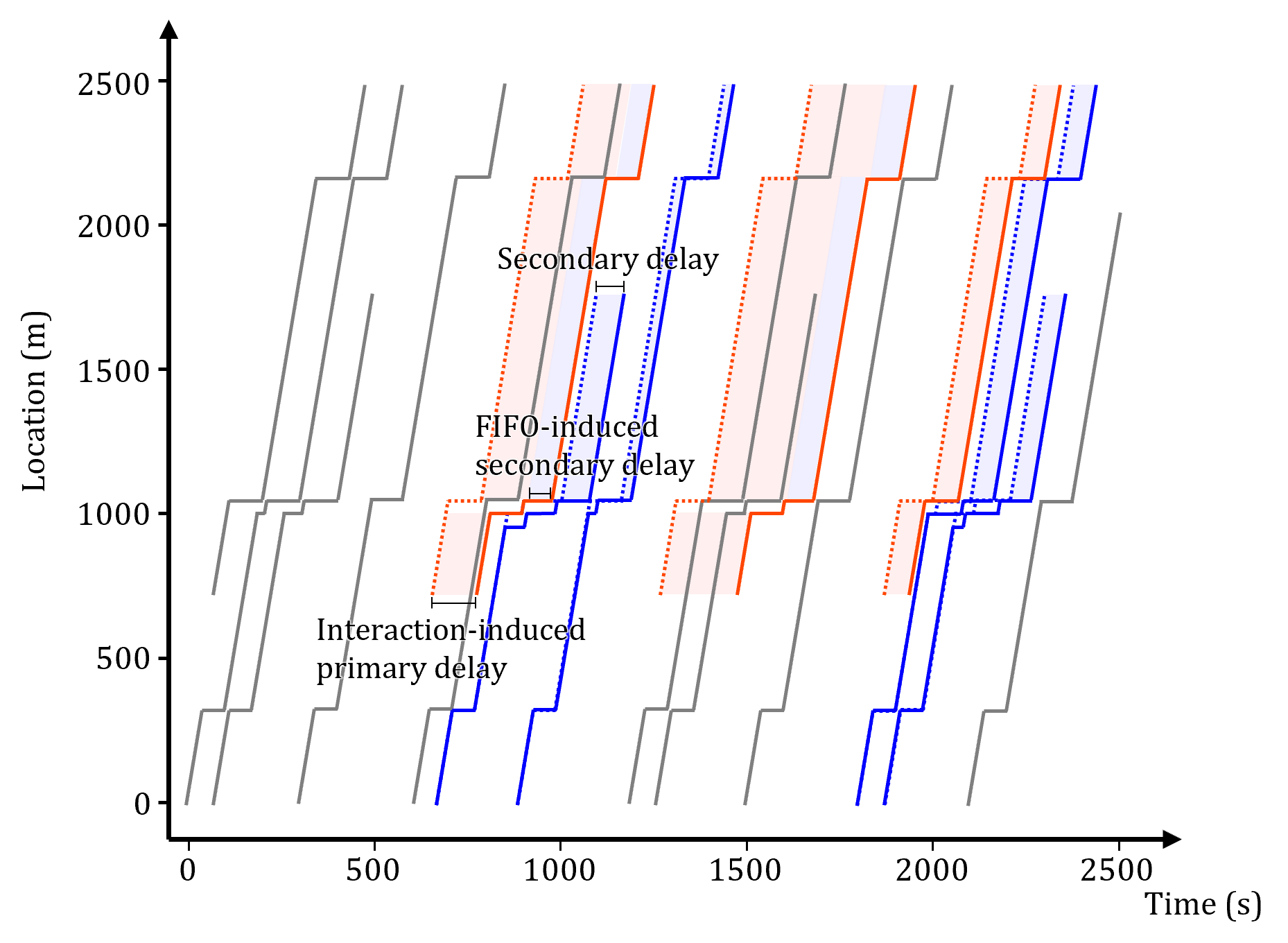}
        \caption{Tram trajectories.}
        
    \end{subfigure}
    \begin{subfigure}[b]{0.49\linewidth}
        \includegraphics[width=\linewidth]{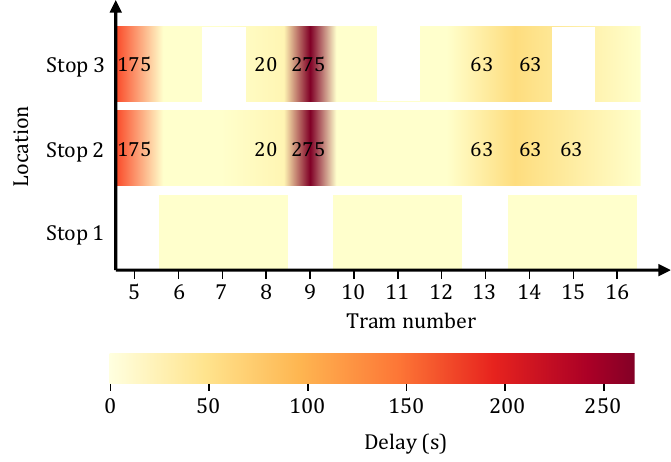}
        \caption{Tramway delay propagation.}
        
    \end{subfigure}
    \caption{Tram trajectories and tramway delay propagation on the eastbound horizontal artery without disturbance. The delayed, influenced by secondary delay, and uninfluenced tram trajectories are distinguished by red, blue, and gray lines. The solid and dashed lines denote the actual and scheduled trajectories.}
    \label{artery_tramway_nodelay}
\end{figure}

Figure \ref{artery_tramway_nodelay} (a) depicts the tram trajectories in the first 2500 s along the eastbound horizontal artery derived from the tramway cumulative flow in M-LTM, assuming trams operate at a constant speed on dedicated tracks. Although there is no disturbance of trams before they enter the network, the interactions with road traffic on the vertical link and signal control cause congestion on the north link, which results in delayed arrivals of trams of Line 2 after the first 500 s. The delay is amplified and leads to secondary delays on the subsequent trams according to the FIFO rule. The tramway delay propagation on partial trams is demonstrated by the heatmap in Figure \ref{artery_tramway_nodelay} (b). Along the horizontal artery, the delays of individual trams increase at Stop 2, while the delays remain constant between Stop 2 and Stop 3. That is because the dedicated ROW of tramway traffic guarantees the stable operation, but Line 2 merges with the other two trams so that more trams meet at Stop 2.

\begin{figure}[!h]
    \centering
    \begin{subfigure}[b]{0.49\linewidth}
        \includegraphics[width=\linewidth]{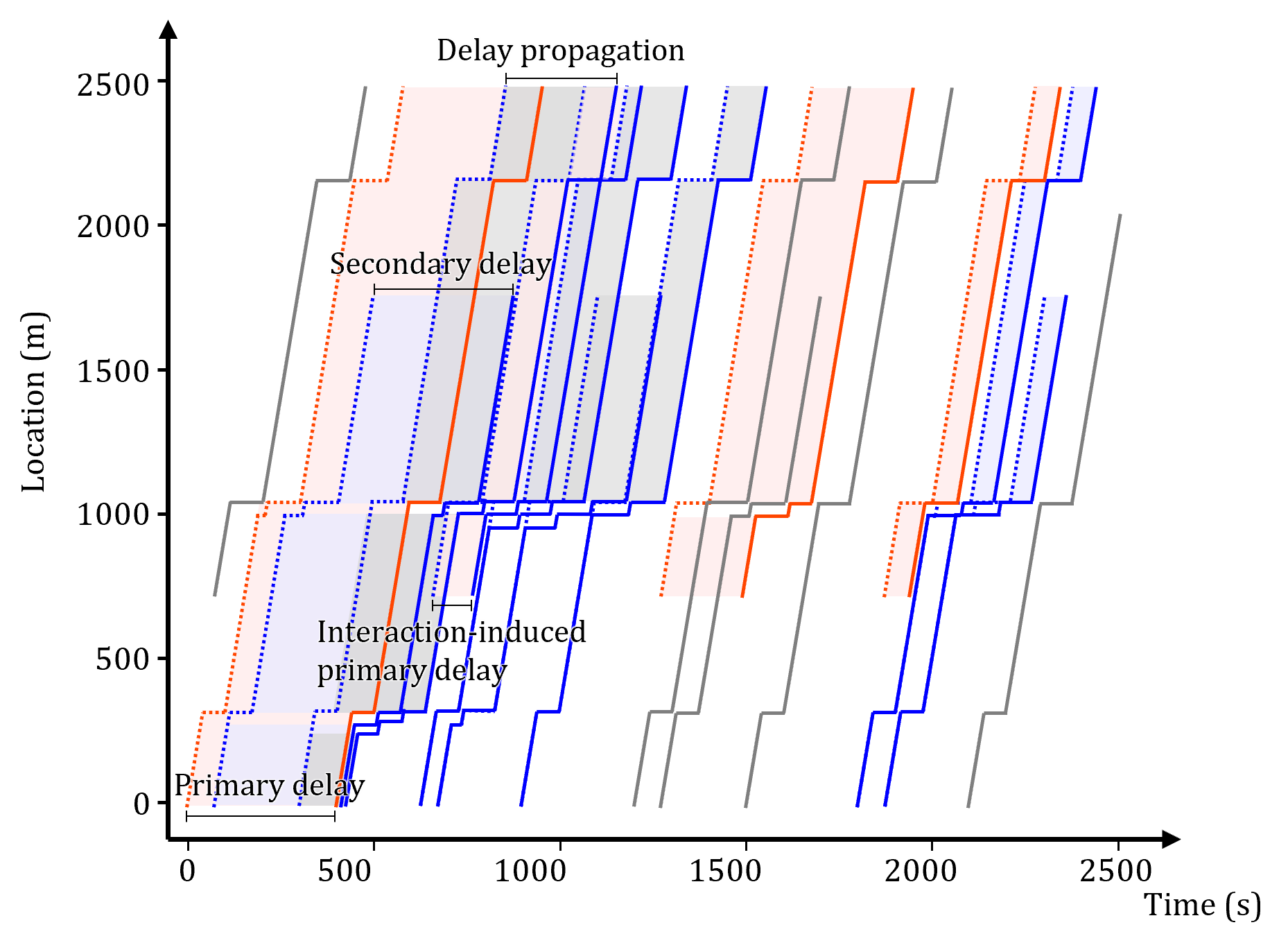}
        \caption{Tram trajectories.}
        
    \end{subfigure}
    \begin{subfigure}[b]{0.49\linewidth}
        \includegraphics[width=\linewidth]{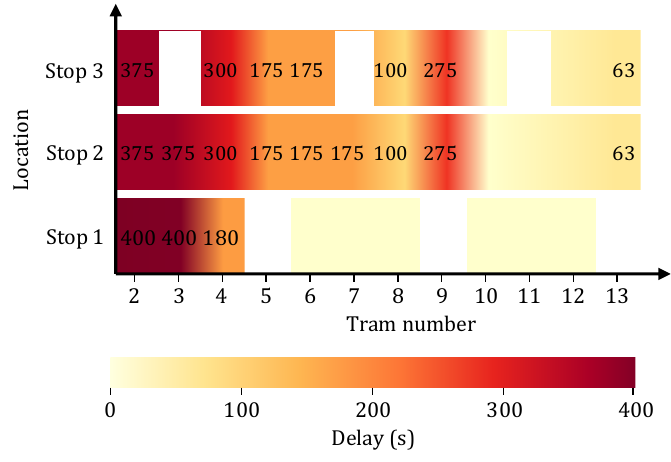}
        \caption{Tramway delay propagation.}
        
    \end{subfigure}
    \caption{Tram trajectories and tramway delay propagation on the eastbound horizontal artery with a disturbance of 400 s. The delayed, influenced by secondary delay, and uninfluenced tram trajectories are distinguished by red, blue, and gray lines. The solid and dashed lines denote the actual and scheduled trajectories.}
    \label{artery_tramway_delay400}
\end{figure}

\begin{figure}[!h]
    \centering
    \begin{subfigure}[b]{0.49\linewidth}
        \includegraphics[width=\linewidth]{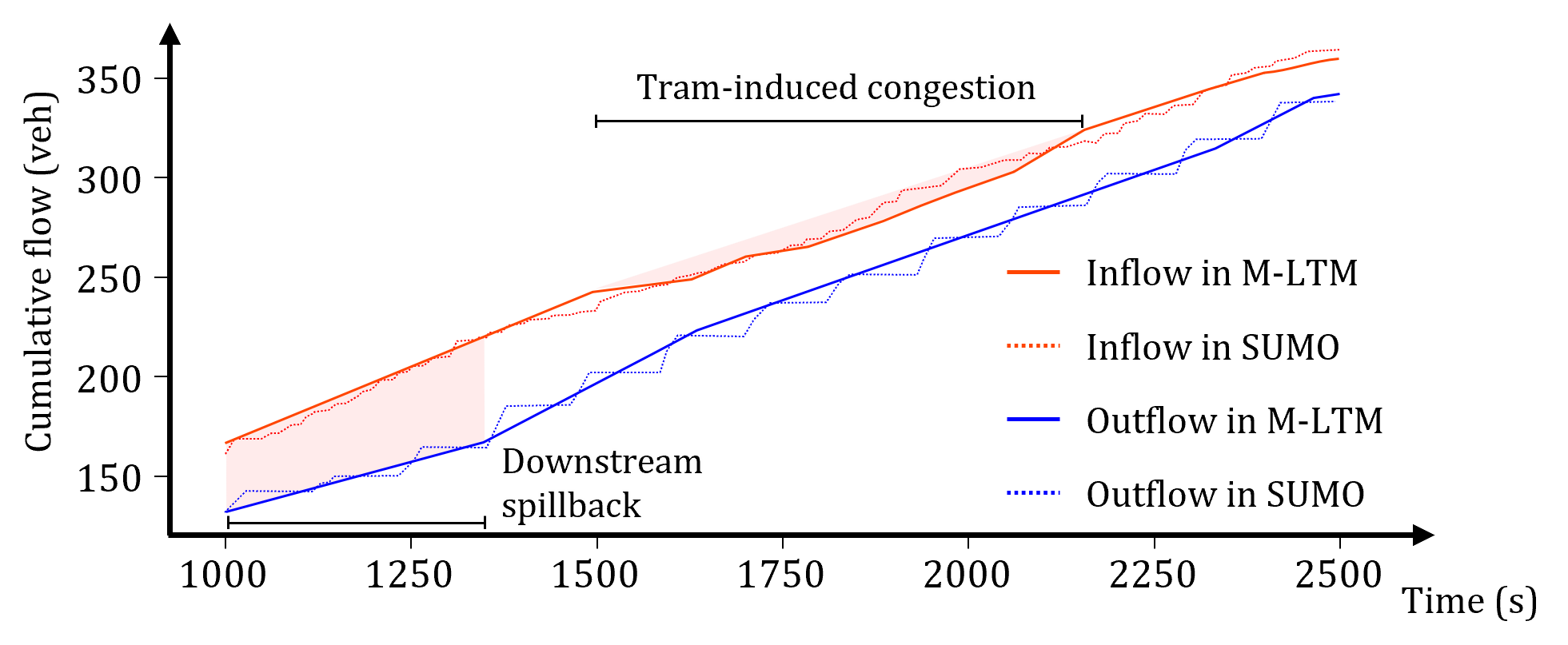}
        \caption{Southbound vertical link.}
        
    \end{subfigure}
    \begin{subfigure}[b]{0.49\linewidth}
        \includegraphics[width=\linewidth]{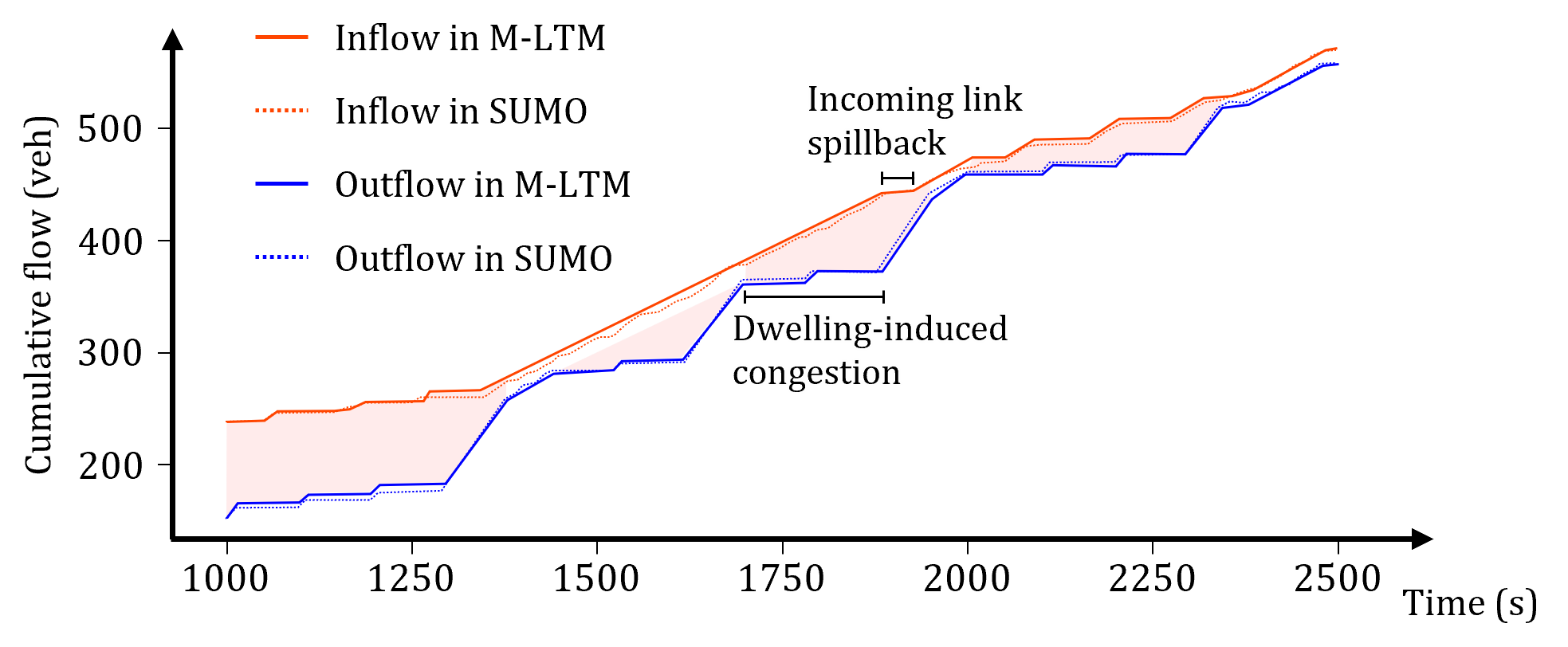}
        \caption{Eastbound incoming link of Stop 2.}
        
    \end{subfigure}
    \caption{Road traffic cumulative flow curve on the southbound vertical link and the eastbound incoming link of Stop 2 in the scenario of no tram disturbance.}
    \label{road traffic cumulative flow arterial}
\end{figure}

A 400 s disturbance is assumed on the entry of the first tram of Line 1 to further illustrate the ability of the M-LTM to analyze tramway delay evolution, as shown by the trajectories and delay heatmap in Figure \ref{artery_tramway_delay400}. The delay first propagates to the entry of the subsequent trams, while the primary delay reduces later. There is no tram at Stop 2 when the delayed tram arrives, avoiding additional waiting before the stop as scheduled. However, it further induces secondary delays at the stop because of the FIFO rule. Additionally, the interaction-induced delay of Line 2 aggregates the congestion at the stop. Compared to the road traffic congestion-induced delay, the disturbance results in a severe negative impact on the following trams. This suggests that the disturbances of trams may cause the multimodal traffic system to be unstable, highlighting the necessity of coordinated management and control under disturbances.
The road traffic results are compared with SUMO simulation and evaluated by the normalized deviation of cumulative flow at the upstream and downstream ends of links and the GEH statistic of each link. The average normalized deviations in the two scenarios are 5.93\% and 9.82\%, and the maximum GEH statistic is 3.65. To show the detailed results, the road traffic cumulative flow during the 1000 s and 2500 s of the southbound vertical link and the eastbound incoming link of Stop 2 are exhibited in Figure \ref{road traffic cumulative flow arterial}. The tram dwellings result in road traffic congestion on the stop incoming link, which propagates backward and leads to link spillback. The spillback further impedes the inflow of the outbound vertical link, increasing the road traffic density. Therefore, the trams reduce the link capacity as moving bottlenecks.

\subsection{Dresden network}
\label{dresden_case}
The proposed M-LTM is applied to the center of the Dresden network to illustrate the multimodal traffic dynamics and the application of the model in complex urban networks.

\subsubsection{Model verification on Dresden network}
\label{dresden_verificiation}
The Dresden network used in this case study consists of 100 links, 24 signalized intersections, and 70 tram stops. A total of 11 tram lines operate in the study area, covering 68 links. Among these, 54 links are equipped with dedicated tram tracks, while the remaining 14 links operate under shared ROW conditions for road traffic and tramway traffic.

Road traffic demands are injected from 15 source nodes based on the real-world data during morning peak hours (from 6:00 AM to 8:00 AM) on Monday, 8th September 2025, obtained from the VAMOS system. There are no external disturbances in tramway traffic and no remarkable incidents in the road traffic system on the day. The turning ratio, saturation flow at intersections, and free flow speed of each link are calibrated based on network layout and the real-world data on weekdays from 1st to 14th September, except for the simulated day. The VAMOS system is a traffic management system in Dresden, providing real-time traffic measurements through inductive loops, Traffic Eye Universe (TEU) cameras, parking detectors, and roadside units (RSUs) \citep{ijaradar2026did}. Besides the detectors at the source nodes, there are 46 detectors located in the research area. The traffic flow data collected from these detectors is used to validate the road traffic dynamics generated by the proposed model. The tram position data collected by on-board devices on vehicles provides information for the tramway traffic demands from source nodes, and provides ground truth for validating the simulated tramway traffic dynamics. The tram dwelling times are obtained from tram timetables.

\begin{figure}[!h]
    \centering
    \begin{subfigure}[b]{0.49\linewidth}
        \includegraphics[width=\linewidth]{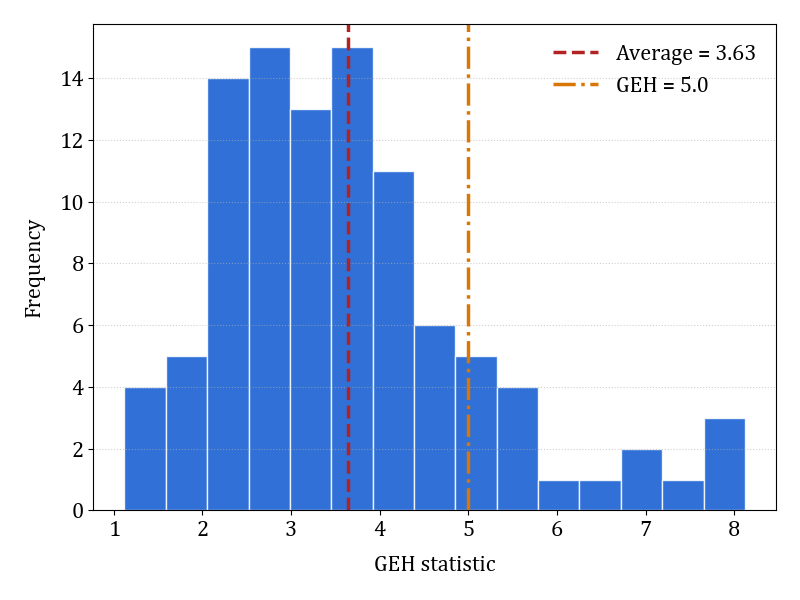}
        \caption{The frequency distribution of GEH statistic of road traffic across the Dresden network.}
        
    \end{subfigure}
    \begin{subfigure}[b]{0.49\linewidth}
        \includegraphics[width=\linewidth]{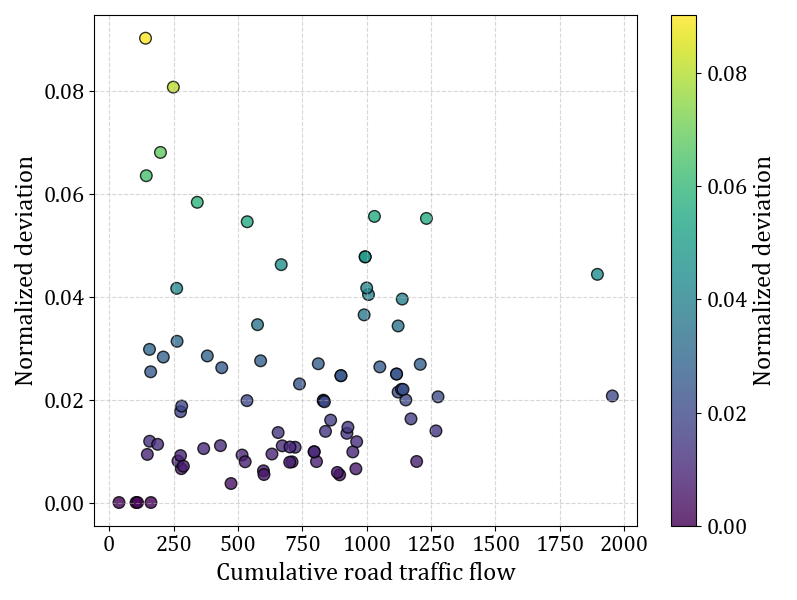}
        \caption{The normalized deviation distribution with cumulative road traffic outflow on the Dresden network.}
        
    \end{subfigure}
    \caption{GEH statistic and normalized deviation distribution of road traffic cumulative flow on the Dresden network.}
    \label{Dresden_normalized_deviation}
\end{figure}

The accuracy of the proposed model is evaluated by comparing the cumulative road traffic flow with real-world data from VAMOS system at detectors. The average GEH statistic and normalized deviation of cumulative road traffic flows over all detectors are 3.63 and 2.30\%, suggesting that the proposed model accurately captures the road traffic dynamics across the whole network. The frequency distribution of the GEH is depicted in Figure \ref{Dresden_normalized_deviation} (a). 86\% road traffic flows have a GEH less than 5.0. The relationship between the normalized deviation and cumulative road traffic flow is depicted in Figure \ref{Dresden_normalized_deviation} (b). Relatively large deviations tend to occur at locations with low cumulative traffic flows. This phenomenon can be attributed to two factors. First, traffic flow is inherently more heterogeneous under low-demand conditions, which makes accurate reproduction more challenging. Second, links with small measured flows often play minor roles in the arterial structure and are more strongly influenced by local streets and unmodeled access flows, which are not explicitly represented in this case study.

The tram arrival and departure times at stops derived from the proposed model are compared with the real-world data. The average absolute errors between simulation results and real-world data of arrival times and departure times are 8.73 s and 9.46 s, respectively. These small discrepancies demonstrate that the tram dwelling process and the temporal evolution of tram operations are accurately captured by the proposed model.

\begin{figure}[!h]
   \begin{subfigure}[b]{0.45\linewidth}
        \includegraphics[width=\linewidth]{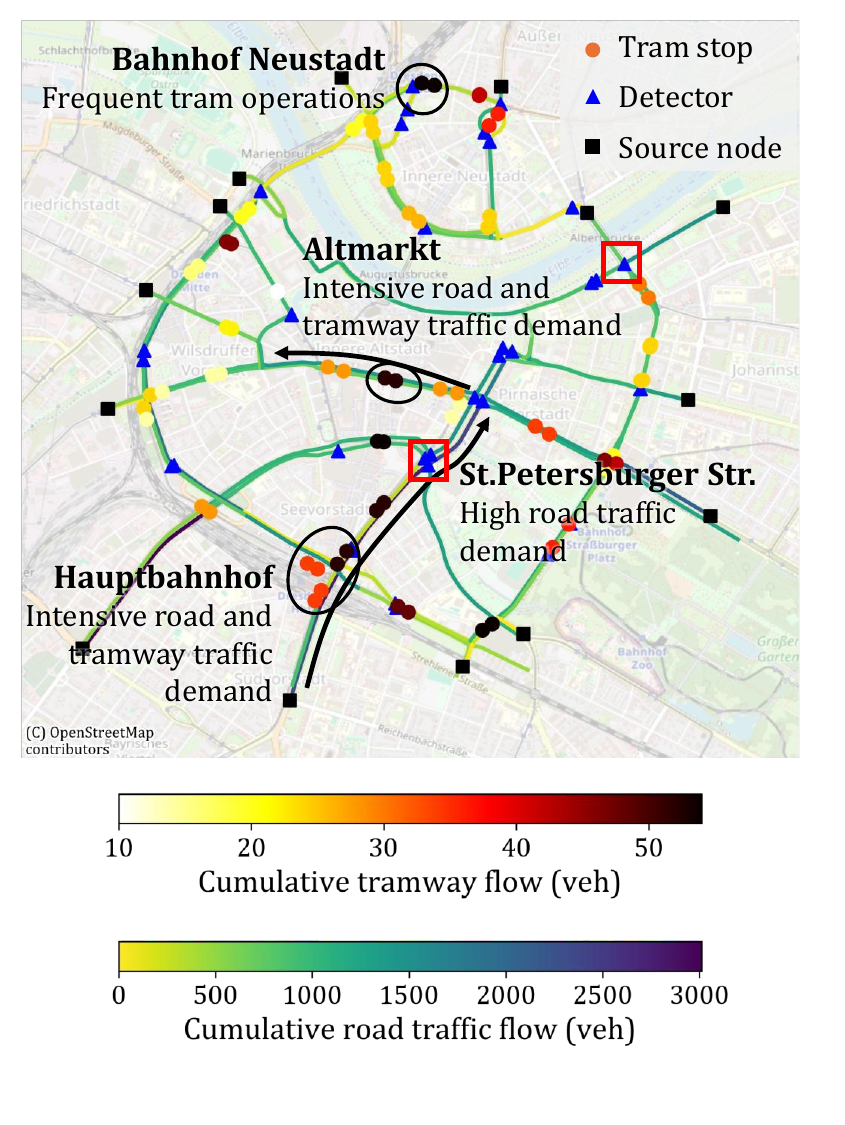}
        \caption{Cumulative flow.}
        \label{cumulative_flow_dresden}
    \end{subfigure}
     \begin{subfigure}[b]{0.45\linewidth}
        \includegraphics[width=\linewidth]{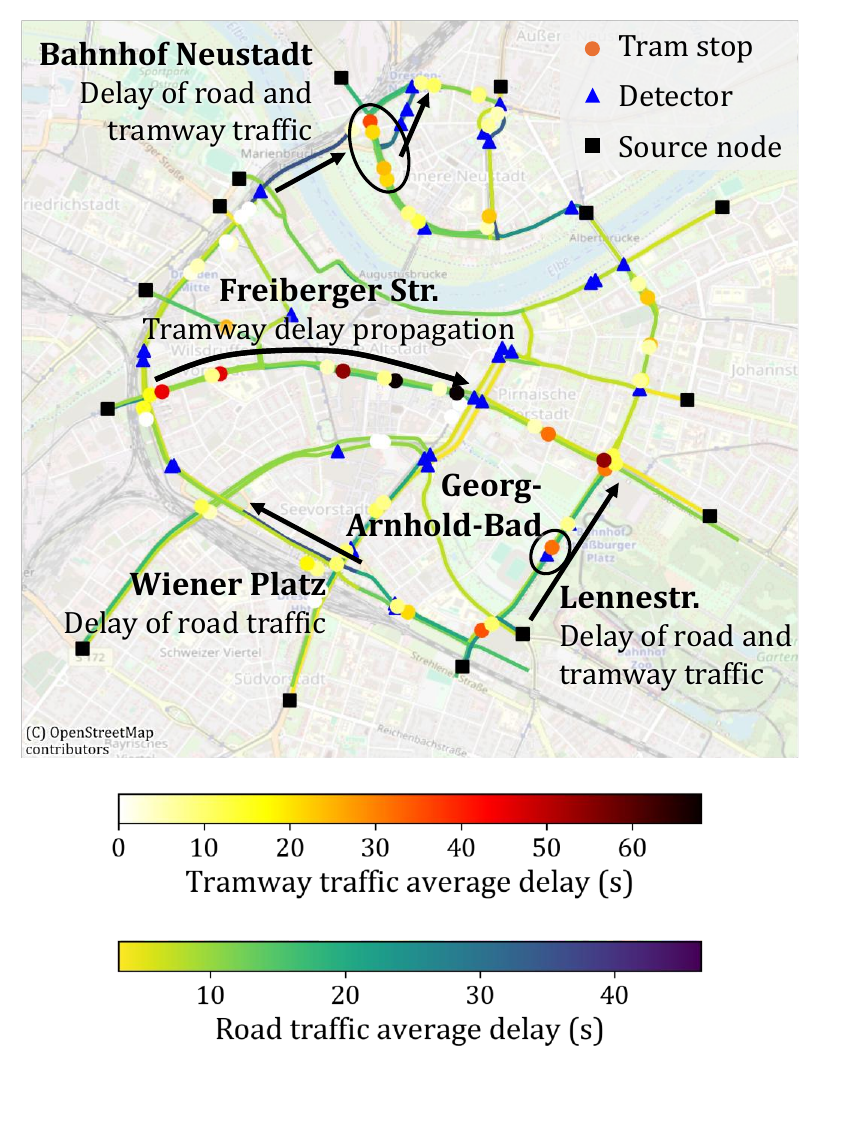}
        \caption{Average delay.}
        \label{average_delay_dresden}
    \end{subfigure}
    \caption{Heatmaps of cumulative traffic flow and average delay of road traffic and tramway traffic on the Dresden network.}
    \label{Dresden_heatmap}
\end{figure}

\begin{figure}[!h]
   \begin{subfigure}[b]{0.45\linewidth}
        \includegraphics[width=\linewidth]{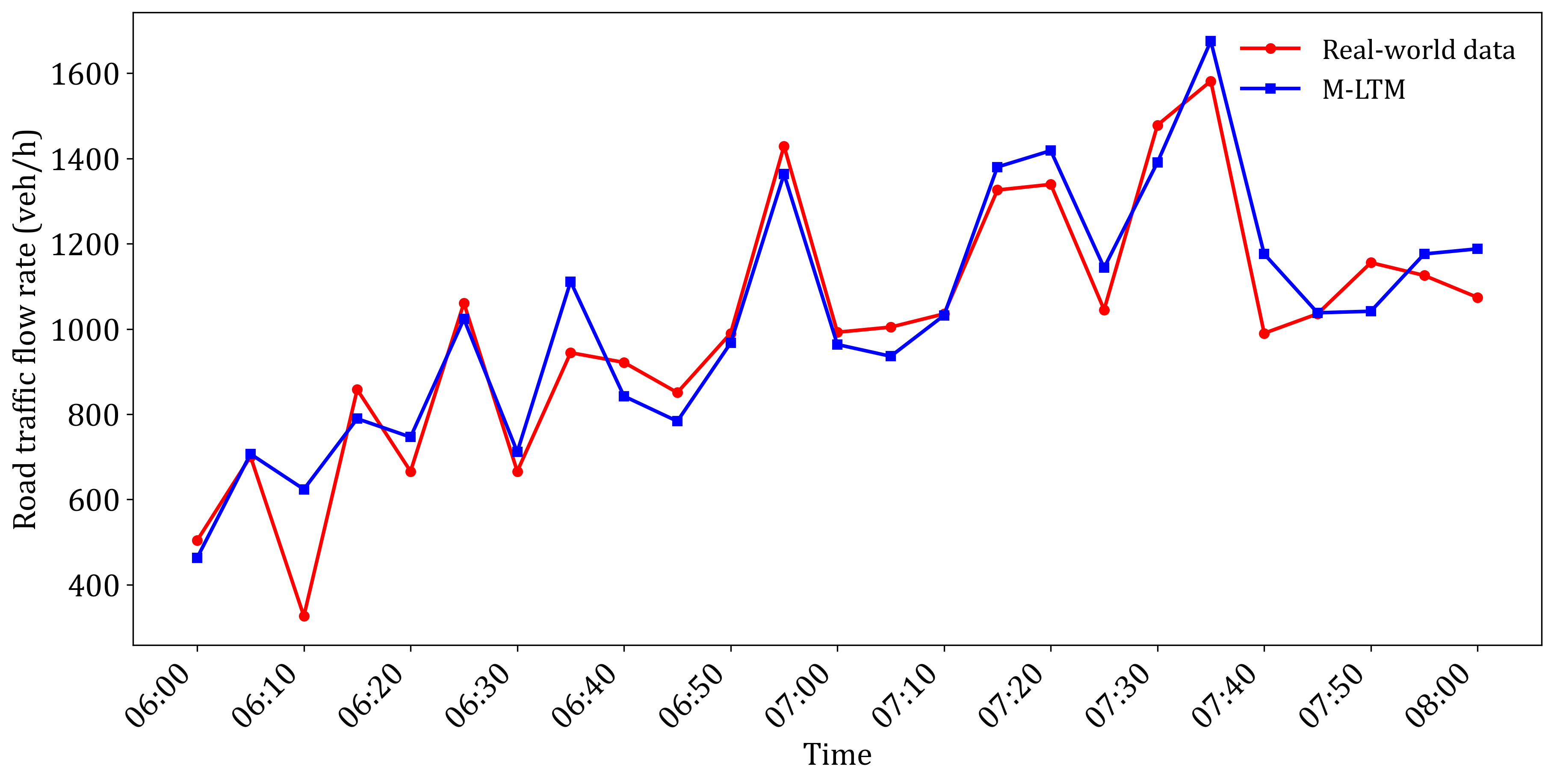}
        \caption{Road traffic flow rate on St.Petersburger Str.}
        \label{Road traffic flow rate on St.Petersburger Str.}
    \end{subfigure}
     \begin{subfigure}[b]{0.45\linewidth}
        \includegraphics[width=\linewidth]{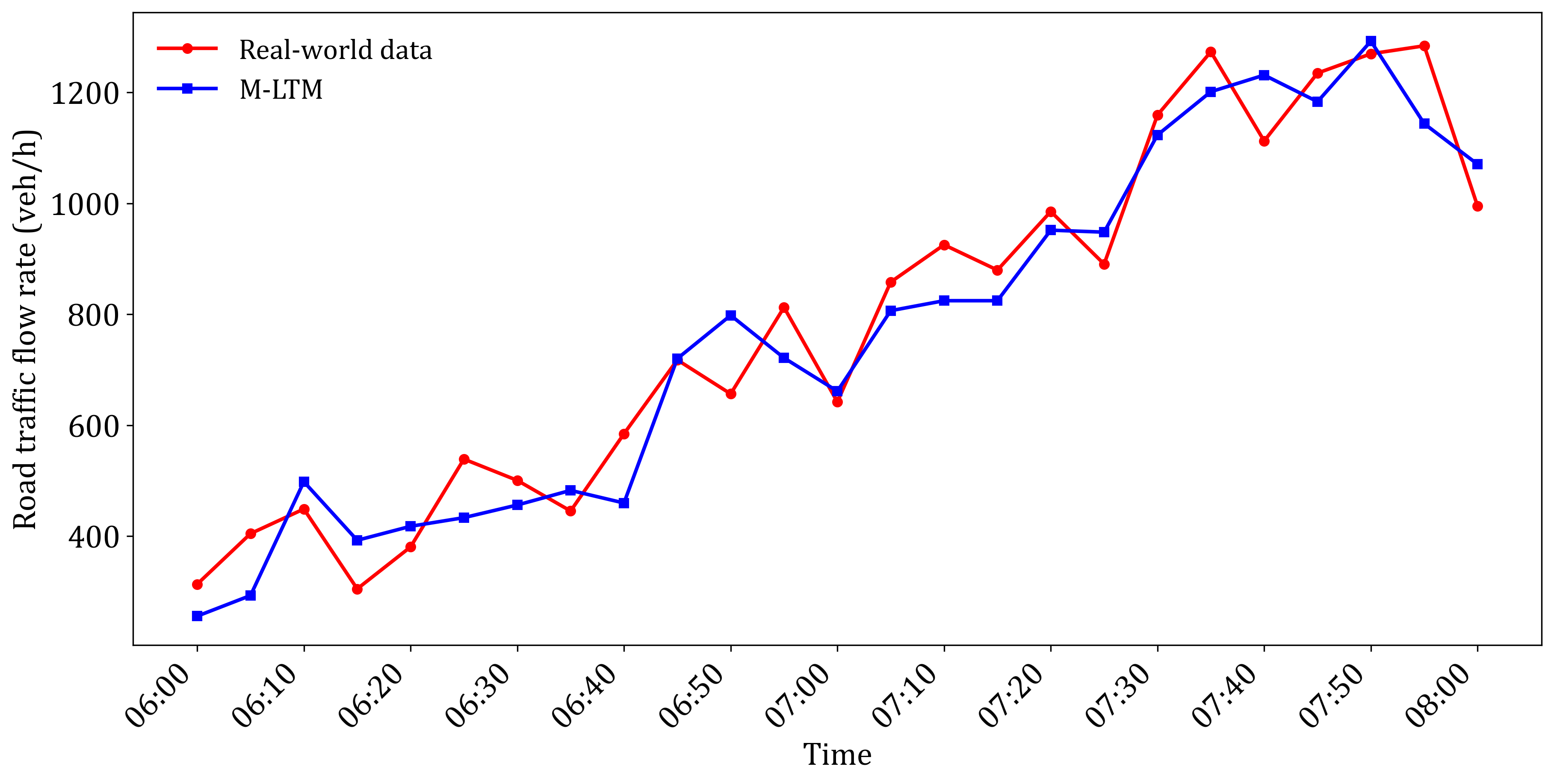}
        \caption{Road traffic flow rate on Albetbrücke.}
        \label{Road traffic flow rate on St.Petersburger Str.}
    \end{subfigure}
    \caption{Road traffic flow rate evolution comparison on St.Petersburger Str. and the Albetbrücke.}
    \label{Dresden_flow_rate}
\end{figure}

To further explore the traffic evolution on the network, we demonstrate the cumulative flow and average delay distribution across the network in Figure \ref{Dresden_heatmap}. In the figures, link colors represent road traffic quantities, while the colors of points at tram stops indicate tramway traffic quantities. Road traffic average delay is calculated as the difference between actual travel time and the minimum travel time at free flow speed, whereas tramway delay is computed as the deviation between scheduled and actual departure times from stops. St. Petersburger Str. is identified as one of the most heavily road traffic-loaded corridors, with cumulative flow reaching up to 2500 vehicles during the two-hour peak period. And Bahnhof Neustadt exhibits the highest tramway demand, with approximately 50 trams passing through. Moreover, two major multimodal hotspots, Hauptbahnhof and Altmarkt, are characterized by simultaneously high volumes of both road traffic and tramway traffic. To further evaluate the accuracy and ability to reproduce traffic dynamics, the road traffic flow rate evolution on St.Petersburger Str. and the upstream of Albetbrücke, whose locations are denoted by red boxes in Figure \ref{Dresden_heatmap} (a), is compared in Figure \ref{Dresden_flow_rate}.  St.Petersburger Str. is located in the middle of the network, with high road traffic demand and isolated tram operations. The consistent flow rate evolution of real-world data and M-LTM demonstrates the model's ability to capture heavy road traffic dynamics. Although tramway traffic has dedicated tracks along the upstream link of Albetbrücke in the northeast of the network, it interacts with road traffic at upstream nodes. The road traffic flow rate from M-LTM has the same trend as the real-world data, indicating that the model can accurately capture the inter-modal interactions.

The spatial distribution of average delay reveals more complex multimodal interactions. Road traffic delay is strongly correlated with traffic flow, whereas tram delay varies significantly across links and is mainly induced by congestion on shared ROWs. Dedicated ROWs effectively maintain tram punctuality, as evidenced near Wiener Platz and Haptbahnhof. Conversely, congestion propagation from road traffic to tram operations is evident along Freiberger Str.. The impacts of multimodal traffic interactions are amplified at the stop Georg-Arnhold-Bad on northbound Lennestr., where the two traffic modes share the ROW over a distance of 110 m between the stop and the downstream intersection. The trams act as moving bottlenecks and render the road traffic unable to change lanes, leading to additional road traffic delay. Contrarily, road traffic waiting at the intersection leads to tramway delays. Similar bidirectional interactions are observed near Bahnhof Neustadt. These indicate the interactions between road and tramway traffic in the multimodal traffic system.
\subsubsection{Model applications on Dresden network}

We further explore the application of the proposed M-LTM in predicting macroscopic control strategy performance. The Park and Ride (P+R) and urban speed limit strategies are implemented, and the resulted traffic changes are derived from the M-LTM.

\begin{figure}[!h]
   \begin{subfigure}[b]{0.45\linewidth}
        \includegraphics[width=\linewidth]{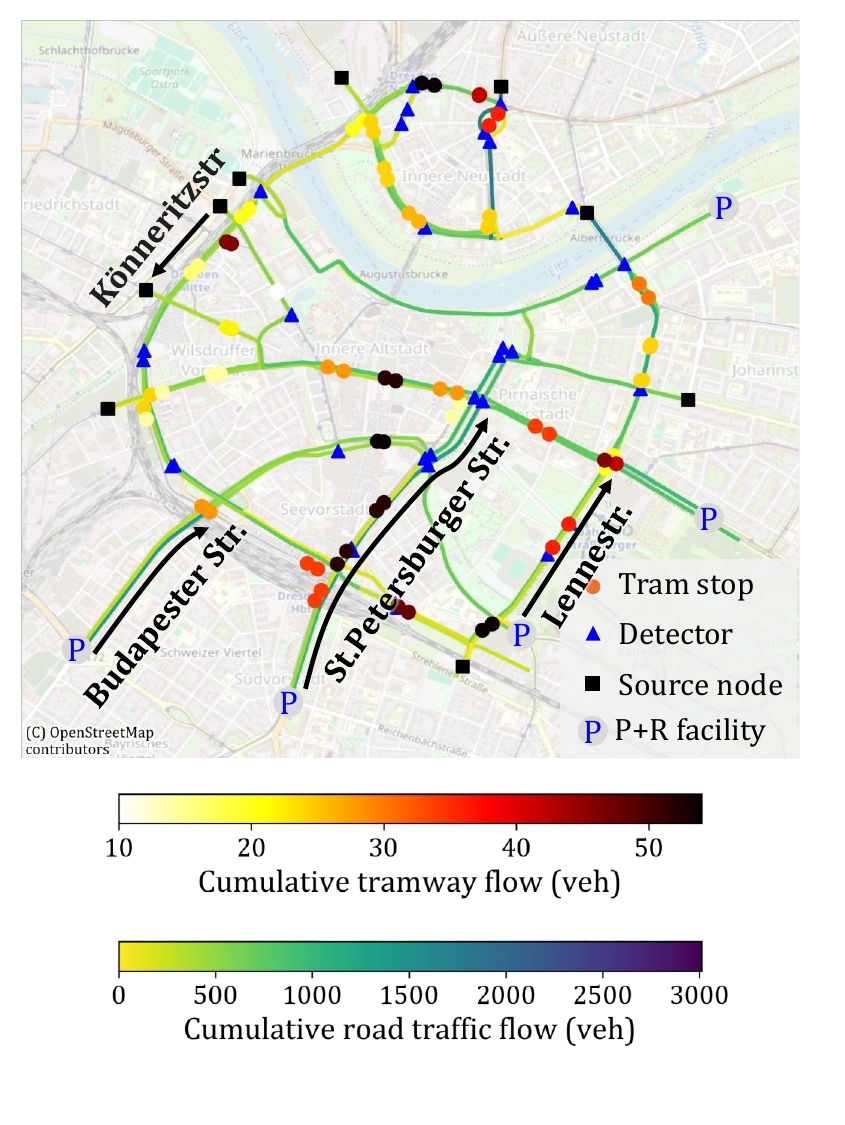}
        \caption{Cumulative flow with P+R strategy.}
        \label{cumulative_flow_pr}
    \end{subfigure}
     \begin{subfigure}[b]{0.45\linewidth}
        \includegraphics[width=\linewidth]{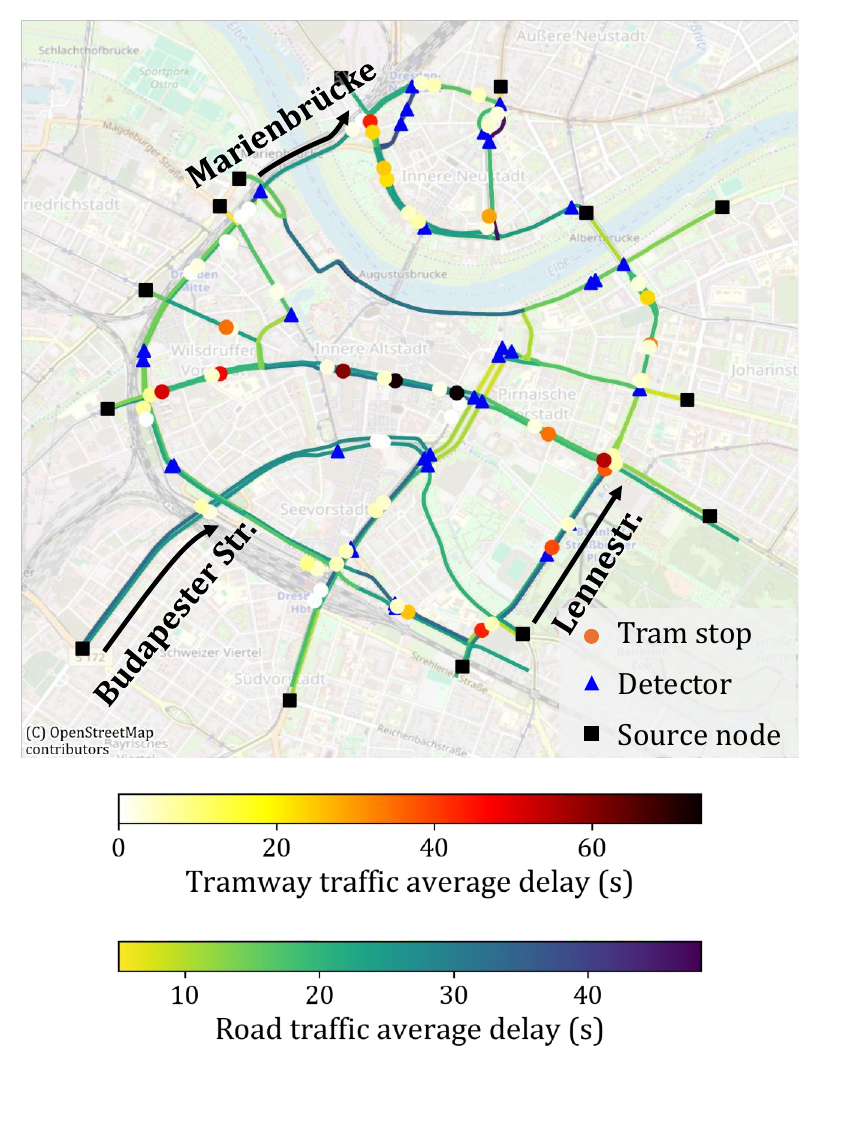}
        \caption{Average delay with urban speed limit strategy.}
        \label{average_delay_speed_limit}
    \end{subfigure}
    \caption{Heatmaps of cumulative traffic flow after P+R strategy and average delay of road traffic and tramway traffic after urban speed limit strategy on the Dresden network.}
    \label{Dresden_control_heatmap}
\end{figure}

\begin{figure}[!h]
    \centering
    \includegraphics[width=0.5\linewidth]{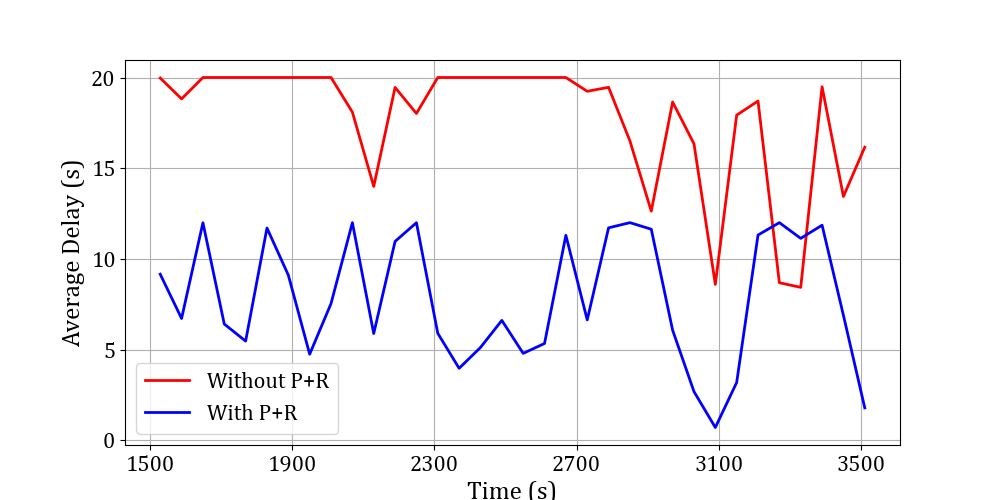}
    \caption{Road traffic delay evolution during 1500 s (6:25 AM) and 3600 s (7:00 AM) on northbound Lennestr..}
    \label{delay_pr}
\end{figure}

The P+R strategy is implemented at five source nodes, as denoted in Figure \ref{cumulative_flow_pr}, to reduce road traffic flows in the city center. The road traffic demand at these five source nodes with P+R facilities reduces 80\%, while other road traffic demand and tramway demand remain the same as the original data. A significant redistribution of road traffic is observed after the P+R strategy, particularly with reductions in road traffic flows on Lennestr., St.Petersburger Str., and Budapester Str., as shown in Figure \ref{cumulative_flow_pr}. Focusing on the critical bottleneck with inter-modal interactions at the stop Georg-Arnhold-Bad on northbound Lennestr., the congestion is alleviated notably with the P+R strategy, with reductions of 43.6\% and 56.7\% on road and tramway traffic. Figure \ref{delay_pr} exhibits the road traffic delay evolution during 1500 s (6:25 AM) and 3600 s (7:00 AM) on northbound Lennestr.. Besides the overall reduction, the P+R strategy shortens the congestion duration. On the one hand, less road traffic is blocked by the tram on the shared-ROW link, and avoid spillback from the incoming link. On the other hand, the tram will not be delayed by the left-turning road traffic at the intersection, which reduces the propagation of congestion on both traffic modes.

Lower urban speed limit has been widely advocated in Europe to enhance traffic safety \citep{lu2023simulation}. The average free flow speed calibrated by the real data is 51.3 km/h, and a new speed of 30 km/h is set in the model. While the lower speed limit leads to more congested multimodal traffic, the road traffic delay distribution tends to be homogeneous across the network, as shown in Figure \ref{average_delay_speed_limit}. The road traffic average delay on the critical link, Marienbrücke, decreases after the control strategy. Conversely, the northbound Budapester Str., which is noticed with slight delay, becomes congested as well. The heavy road congestion aggregates the tramway delay induced by the inter-modal interactions. The tramway operations at stops with marginal delay are not influenced significantly, while delays on shared-ROW links are amplified. Particularly, both the road and tramway traffic delays on the northbound Lennestr. increase after implementing the lower speed limit.

The Dresden case study demonstrates that the proposed M-LTM not only accurately captures the multimodal traffic dynamics but also offers a powerful analytical tool for diagnosing congestion mechanisms and inter-modal interactions in complex urban networks. The ability to detect road traffic congestion, tramway delay, and their interaction processes the potential of the proposed framework for multimodal traffic management, planning, and coordinated control applications.

\section{Conclusion}
\label{conclusion}
This study develops a multimodal link transmission model for road and tramway traffic that captures multimodal traffic dynamics in a unified continuous-time framework. Specifically, the link model describes the traffic propagation under  inter-modal interactions on links with shared ROW. The node model achieves traffic transfer and tram operations at tram stops and signalized intersections. The tram stop node model simulates the tram dwelling process and dwelling-induced congestion of road traffic, satisfying the FIFO rule, while the signalized intersection node distributes multimodal traffic flows according to the invariance and holding-free principles. Moreover, the road traffic spillback at intersections is explicitly considered, further capturing the inter-modal interactions.

The proposed model is implemented on two one-node networks and an arterial network to illustrate detailed multimodal traffic evolutions and interaction mechanisms. Furthermore, a case study on the Dresden network based on real-world data is conducted to validate the model’s accuracy with an average normalized deviation of cumulative flow of only 2.3\% and to demonstrate the ability of M-LTM in traffic analysis and control applications.

The Multimodal Link Transmission Model demonstrates significant advantages in macroscopically simulating the urban traffic dynamics of road traffic and tramway traffic and their interactions. The stochastic factors, such as fluctuations in tram dwelling time and free flow speed, will be addressed in the future. Further research will concentrate on designing coordinated control strategies based on the proposed model. Additionally, incorporating other traffic modes, such as bus and bicycle traffic, into the proposed model will also be explored.

\section*{Acknowledgement}
The data of this study is obtained from the VAMOS system, funded by the federal government and the Free State of Saxony, the traffic management system for the state capital of Dresden, and Dresdner Verkehrsbetriebe AG. This study is also supported by the Graduate Academy of Technische Universität Dresden through Maria Reiche Doctoral Fellowship.

\appendix
\section{Algorithm for determining multimodal traffic internal boundary conditions with interactions}
\label{internal_algorithm_appendix}
\begin{figure}[!h]
    \centering
    \includegraphics[width=\linewidth]{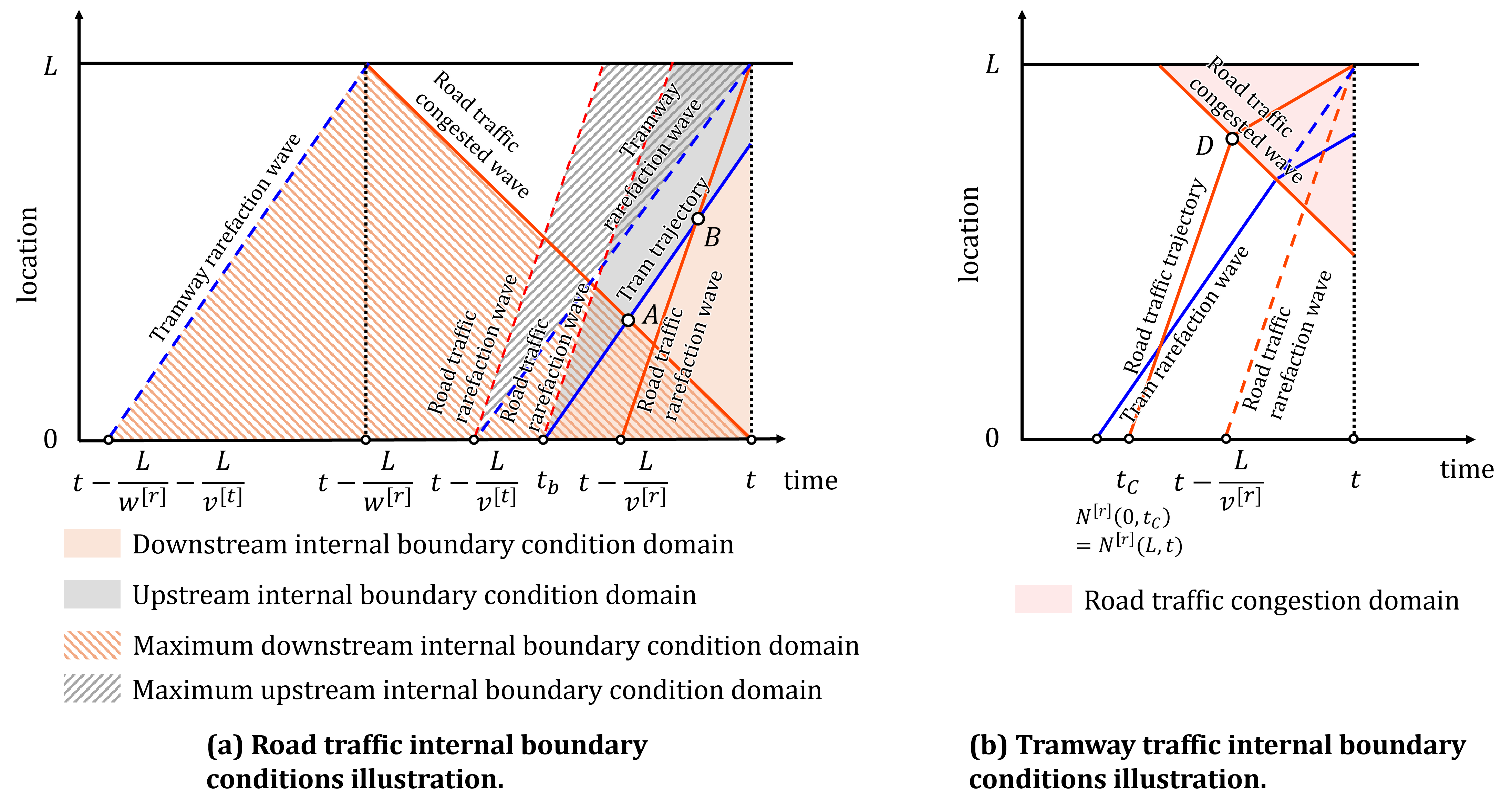}
    \caption{Illustration of internal boundary conditions.}
    \label{road traffic internal boundary illustration}
\end{figure}

When the road and tramway traffic share the ROW on links, they interact with each other, inducing internal boundary conditions in the link model. Figure \ref{road traffic internal boundary illustration} illustrates the traffic wave propagation for the two traffic modes with internal boundary conditions. The processes of determining internal boundary conditions for road traffic and tramway traffic are demonstrated in Algorithm \ref{road traffic internal boundary} and Algorithm \ref{tramway internal boundary condition}, respectively.

\begin{algorithm}[H]
  \caption{Algorithm for determining road traffic internal boundary conditions.}
  \label{road traffic internal boundary}
  \KwIn{Fundamental diagram parameters: $v^{[r]}$, $v^{[t]}$, $w^{[r]}$, $\rho_j^{[r]}$}
  \KwIn{Previous cumulative inflow and outflow: $N^{[r]}(0,t')$, $N^{[t]}(0,t')$, $N^{[r]}(L,t')$, $N^{[t]}(L,t')$, for $t'<t$}
  \KwOut{Road traffic potential cumulative flows imposed by internal boundary conditions: $\bar{N}^{[r]}_1(L,t)$, $\bar{N}^{[r]}_2(0,t)$}
  \textbf{Step 1: Initialize and identify the tram entry time.}\\
  $\bar{N}^{[r]}_1(L,t)=M$, $\bar{N}^{[r]}_2(0,t)=M$, where $M$ is an infinite real number.\\
  $t_b=\max\{t_b|N^{[t]}(0,t_b)<\max \{N^{[t]}(0,t'),t'<t\}\}$\\
  \textbf{Step 2: Determine the road traffic downstream boundary condition.}\\
  \If{$t-\frac{L}{w^{[r]}}-\frac{L}{v^{[t]}}<t_b<t$}{
  $x_A=v^{[t]}w^{[r]}\frac{t-t_b}{v^{[t]}+w^{[r]}}$\\
  $t_A=\frac{v^{[t]}t_b+w^{[r]}t}{v^{[t]}+w^{[r]}}$\\
    $N^{[r]}(x_A,t_A)=\min\{N^{[r]}(0,t-\frac{x_A}{v^{[r]}}),N^{[r]}(0,t_b)+q_r(t_A-t_b),N^{[r]}(L,t-\frac{L-x_A}{w^{[r]}})+\rho_j^{[r]}(L-x_A)\}$\\
  $\bar{N}_d^{[r]}(0,t)=N^{[r]}(x_A,t_A)+x_A\rho_j^{[r]}$\\
  \textbf{Step 3: Determine the road traffic upstream boundary condition.}\\
  \If{$t-\frac{L}{v^{[t]}}\leq t_b< t-\frac{L}{v^{[r]}}$}{
  $x_B=v^{[t]}v^{[r]}\frac{t-t_b}{v^{[r]}-v^{[t]}}$\\
  $t_B=\frac{v^{[t]}t_b-v^{[r]}t }{v^{[t]}-v^{[r]}}$\\
  $N^{[r]}(x_B,t_B)=\min\{N^{[r]}(0,t-\frac{x_B}{v^{[r]}}),N^{[r]}(0,t_b)+q_r(t_B-t_b),N^{[r]}(L,t-\frac{L-x_B}{w^{[r]}}+\rho_j^{[r]}(L-x_B)\}$\\
  \If{$N^{[r]}(x_B,t_B)<N^{[r]}(L,t-\frac{L-x_B}{w^{[r]}}+\rho_j^{[r]}(L-x_B)$}{
  $\bar{N}_u^{[r]}(L,t)=N^{[r]}(x_B,t_B)$
  }
  }
  }
\end{algorithm}

The tramway internal boundary conditions are activated when the surrounding road traffic is congested at a speed lower than that of tramway traffic, calculated by Algorithm \ref{tramway internal boundary condition}. We first discuss whether the road traffic is congested, and identify the congested wave and trajectory of congested road traffic. Then, the tramway traffic potential cumulative flow induced by the upstream internal boundary condition $\bar{N}_u(L,t)$ is derived from the road traffic congested wave. Finally, the restriction of tramway traffic cumulative inflow, $\bar{N}_d(0,t)$, is determined by the link capacity considering road traffic flow.
\begin{algorithm}[H]
\caption{Algorithm for determining tramway traffic internal boundary conditions.}
\label{tramway internal boundary condition}
     \KwIn{Road traffic fundamental diagram: $Q(\rho)$,  $v^{[r]}$, $v^{[t]}$, $w^{[r]}$, $\rho_j^{[r]}$}
  \KwIn{Road traffic cumulative inflow and outflow: $N^{[r]}(0,t')$, $N^{[r]}(L,t')$ $t'\leq t$, and previous tramway traffic cumulative inflow and outflow: $N^{[t]}(0,t')$, $N^{[t]}(L,t')$ $t'<t$}
  \KwOut{Tramway traffic potential cumulative flows imposed by internal boundary conditions: $\bar{N}^{[t]}_1(L,t)$, $\bar{N}^{[t]}_2(0,t)$}
  \textbf{Step 1: Initialize.}\\
  $\bar{N}^{[t]}_1(L,t)=M$, $\bar{N}^{[t]}_2(0,t)=M$\\
  $q^{[r]}(L,t)=\frac{\partial N^{[r]}(L,t)}{\partial t}$\\
  $\rho^{[r]}(L,t)=\max \{Q^{-1}(q^{[r]}(L,t)\}$\\
  $v'=\frac{q^{[r]}(L,t)}{\rho^{[r]}(L,t)}$\\
  \If{$v'<v^{[t]}$}{
   \textbf{Step 2: Determine the tramway internal boundary conditions if the tram is blocked.}\\
  $t_C=\max \{t_C|N^{[r]}(0,t_C)=\max N^{[r]}(L,t'),t'<t\}$\\
  $x_D=\frac{v^{[r]}L-v^{[r]}v'(t-t_C)}{v^{[r]}-v'}$\\
  $t_D=t_C+\frac{L-v'(t-t_C)}{v^{[r]}-v'}$\\
  $\bar{N}^{[t]}_1(L,t)=N^{[t]}(x_D,t_D)$\\
  $N^{[t]}(x_D,t_D)=\min\{N^{[t]}(0,t_D-\frac{x_D}{v^{[t]}}),N^{[t]}(L,t_D-\frac{L-x_D}{w^{[t]}})+\rho_j^{[t]}(L-x_D)\}$\\
  }
  $\bar{N}_d^{[t]}(0,t)=\overleftarrow{N}^{[r]}(t)-\max N^{[r]}(0,t'')+\max N^{[t]}(0,t'), t'<t, t''<t$
\end{algorithm}

\section{Proof of invariant holding-free solution} \label{IHF proof}

Algorithm \ref{flow update without spillback} is an invariant holding-free algorithm defined in \cite{jabari2016node}.
\begin{proof}
    We first prove the holding-free property of the algorithm.

    It is intuitive that the solution is holding-free when all demands are satisfied. Since all traffic flows are demand-constrained, the reduction indices are 1 for all outgoing links. It comes to $q^{[r]}_{ij}(t)=\bar{q}^{[r]}_{ij}(t)$. Therefore $q^{[r]}_i(t)=\bar{q}^{[r]}_i(t)$, which achieves the holding-free property:
    \begin{equation}
        \bar{q}^{[r]}_i(t)-q^{[r]}_i(t)=0
    \end{equation}

    Assuming there are incoming links constrained by downstream supply, it holds:
     \begin{equation}
        q_i^{[r]}=\min_{e^{[r]}_{ij}>0}\{\frac{\bar{q}_j^{[r]}(t)}{\alpha_{ij}}\}<\bar{q}_i ^{[r]}\label{C2}
    \end{equation}
    For the outgoing links of $i$ whose supply constraints are not attained, there is residual supply after the demand assignment of $i$:
    \begin{equation}
        \frac{\overleftarrow{q}_j^{[r]}(t)}{\alpha_{ij}}<\min_{e^{[r]}_{ij'}>0}\{\frac{\overleftarrow{q}_{j'}^{[r]}(t)}{\alpha_{ij'}}\}, \overleftarrow{q}_j^{[r]}(t)>0,e_{ij}>0
    \end{equation}
    Then the receiving flow of the outgoing link $j$ in the next iteration is updated by the residual supply in Step 4 of the Algorithm \ref{flow update without spillback}:
    \begin{equation}
        \overleftarrow{q}_j^{[r]}(t)\leftarrow\overleftarrow{q}_j^{[r]}(t)-q_{ij}^{[r]}(t)=\overleftarrow{q}_j^{[r]}(t)-\alpha_{ij}\min_{e^{[r]}_{ij'}>0}\{\frac{\overleftarrow{q}_{j'}^{[r]}(t)}{\alpha_{ij'}}\}>\overleftarrow{q}_j^{[r]}(t)-\alpha_{ij}\frac{\overleftarrow{q}_j^{[r]}(t)}{\alpha_{ij}}=0
    \end{equation}
Then we focus on the outgoing link $j=\arg \min_{e^{[r]}_{ij'}>0}\{\frac{\bar{q}_{j'}^{[r]}(t)}{\alpha_{ij'}}\}$ of $i$ the minimum is attained, whose turning ratio is greater than 0. If $\overleftarrow{q}_j^{[r]}=0$, the receiving flow of $j$ has been exhausted before assigning traffic demand from $i$ and no additional demand can utilize the link. Otherwise, its receiving flow is updated:
\begin{equation}
        \overleftarrow{q}_j^{[r]}(t)\gets \overleftarrow{q}_j^{[r]}(t)-q^{[r]}_{ij}=\overleftarrow{q}_j^{[r]}(t)-\alpha_{ij}\frac{\overleftarrow{q}_j^{[r]}(t)}{\alpha_{ij}}=0
    \end{equation}
    Under this condition, no supply of $j$ can be utilized after the demand assignment of $i$. Consequently, Algorithm \ref{flow update without spillback} produces solutions that satisfy the supply constraints:
    \begin{equation}
        \exists q,\overleftarrow{q}_j^{[r]}-q_j^{[r]}=0
    \end{equation}
    Therefore, Algorithm \ref{flow update without spillback} produces holding-free solutions.
    
Then, we will prove that the solution is invariant. Suppose the incoming link $i$ is supply constrained, for any demand $\overrightarrow{q}_i^{'[r]}$ greater than $\overrightarrow{q}_i^{[r]}$, it comes:
  \begin{equation}
        \min\{\min_{j\in O_{i}}\{\frac{\overleftarrow{q}_j^{[r]}(t)}{\alpha_{ij}}\},\min\{\overrightarrow{q}_{i}^{'[r]},s_{i}^{[r]}\}\}=\min\{\min_{j\in O_{i}}\{\frac{\overleftarrow{q}_j^{[r]}(t)}{\alpha_{ij}}\},\min\{\bar{q}_{i}^{[r]},s_{i}^{[r]}\}\}=q_{i}^{[r]}
    \end{equation}
    This means that the demand increase will not lead to changes in actual traffic flow. On the other hand, if the traffic flow is demand constrained and there are outgoing links with residual supply after all demand assignment, the increase in their supply has no impact on the solution.

    In conclusion, Algorithm \ref{flow update without spillback} is an invariant holding-free algorithm defined in \cite{jabari2016node}.
        
\end{proof}

\bibliographystyle{elsarticle-harv} 
\bibliography{reference}

\end{document}